
\magnification=1000
\font\tit=cmbx10 scaled\magstep2
\def\P{{\bf P}}
\def\Pn#1{{\bf P}^{#1}}
\def\Pnd#1{{\check{\bf P}}^{#1}}
 

\def\coh#1#2#3{h^{#1}(\O_{#2}{#3})}

\def\Sp{{\bf S}}
\def\Gr{{\bf G}} 
\def\Q{{\bf Q}} 
\def\C{{\bf C}} 
\def\Z{{\bf Z}} 
\def\O{{\cal O}} 
\def\bC{{\bf C}}
\def\bG{{\bf G}}
\def\bP{{\bf P}}

\def\slutt{\sqcap\!\!\!\!\sqcup}

\vskip40pt
{\centerline{ \tit Varieties of sums of powers }}

\vskip10pt
\centerline {\bf Dedicated to the memory of Alf B. Aure and Michael Schneider}

\vskip20pt
\centerline{\it By Kristian Ranestad\footnote *{partially supported by the Norwegian Research
Council} at Oslo  and Frank-Olaf Schreyer at Bayreuth}

\vskip20pt
{\bf Abstract}.  The variety of sums of powers of a homogeneous polynomial of degree $d$ in
$n$ variables is defined and investigated in some examples, old and new.   These varieties are
studied via apolarity and syzygies.
  Classical results (cf. [Sylvester 1851],
[Hilbert 1888], [Dixon, Stuart 1906]) and some more recent results of Mukai (cf. [Mukai
1992]) are presented together with new results for the cases $(n,d)=(3,8), (4,2), (5,3)$. In
the last case the variety of sums of 8 powers of a general cubic form is a Fano 5-fold of
index 1 and degree 660. \medskip \vskip20pt
\vbox{\settabs 2 \columns
 \+ 0. Introduction & \cr 
 \+ 1. Apolarity and syzygies & \cr
 \+ 2. Duality, projections and secants & \cr
 \+ 3. Another example & \cr
 \+ 4. Cubic threefolds, first properties  & \cr
\+ 5. Syzygies of 8 general points in $\Pn 4$ & \cr
 \+ 6. Equations and geometry of the spinor varieties ${\Sp_{ev}}$ and ${\Sp_{odd}}$ & \cr
\+ 7. The apolar Artinian Gorenstein ring of a general cubic threefold & \cr
\+ 8. Proof of the main results & \cr
 \+ 9. Invariants of $VSP(F,8)$ & \cr
  \+ References & \cr }
\vskip20pt

{\bf 0. Introduction}

{\vskip10pt \bf 0.1 } Let $f \in \C[x_0,\ldots,x_n]$ be a homogeneous form of degree d. $f$ can be written as a sum
of powers of linear forms
$$f = l_1^d+ \ldots + l_s^d$$
for $s$ sufficiently large. Indeed, if we identify the map $l \mapsto l^d$ with the $d^{th}$ Veronese
embedding $\Pn n \hookrightarrow \Pn {N_d}$, where $N_d={n+d \choose n}-1$, this amounts to say 
that the image spans $\Pn {N_d}$.  What is the minimal possible value of $s$?

A dimension count shows that
$$ s \ge \lceil {1 \over n+1}{n+d \choose n} \rceil $$
for a {\bf general} $f$. With a few exception for $(d,n)$ equality holds by a result of Alexander
and Hirschowitz [ - 1995] and Terracini's Lemma [ - 1911]:
{\vskip10pt \proclaim 0.2 Theorem. {\rm (Alexander-Hirschowitz)}  A general form $f$ of degree $d$ in $n+1$ variables
is a sum of $\lceil {1 \over n+1}{n+d \choose n} \rceil$ powers of linear forms, unless\smallskip 
 $d=2$, where $s=n+1$ instead of $\lceil {n+2 \over 2}\rceil$, or\smallskip
 $d=4$ and $n=2,3,4$, where $s=6,10,15$ instead of $5,9,14$ respectively, or\smallskip
$d=3$ and $n=4$, where $s=8$ instead of $7$.\par

\bigskip 
The exceptions were classically known cf.  [Clebsch 1861], [Reye 1874], 
[Sylvester 1886], [Richmond 1902], [Palatini 1903], [Dixon 1906].  See also 1.4.

{\vskip10pt \bf 0.3 } For {\bf special} forms $f$ the minimal value $s$ can both be smaller or bigger (cf. 1.6). As 
far as we know it is an open problem to determine the function $\tilde s = \tilde s (d,n)$ such that
{\bf every} form of degree d in n+1 variables can be written as a sum of less than $\tilde s+1$
powers. 

{\vskip10pt \bf 0.4 } Given the answer to the first question, one might ask whether the presentation as a sum of
powers is unique. Of course a {\bf general} form can have a unique expression 
only if ${1 \over n+1}{n+d \choose n}$ is an integer. The following is known:
\bigskip

 {\it A general form $f$ of degree $d$ in $n+1$ variables
has a unique presentation as a sum of $s={1 \over n+1}{n+d \choose n} $ powers of linear forms
in the following cases:\smallskip
$n=1$, $d=2k-1$ and $s=k$, {\rm [Sylvester 1851]}, or\smallskip
$n=3$, $d=3$ and $s=5$ {\rm Sylvester's Pentahedral Theorem [Sylvester 1851] [Reye 1873]},
or\smallskip
$n=2$, $d=5$ and $s=7$ {\rm [Hilbert 1888], [Richmond 1902], [Palatini 1903]}.
}
\bigskip
{\bf Question:}  Are there further examples?  \par
Iarrobino and Kanev treat this question for polynomials with fewer summands (cf. [Iarrobino, Kanev
1996]).

{\vskip10pt \bf 0.5 } We put this question in a somewhat more general frame work:
Let $F =V(f) \subset \Pn n$ be a hypersurface of degree $d$. If the equation
$f = l_1^d+\ldots+l_s^d$ is a sum of s powers of linear forms, then the linear forms give
hyperplanes in $ L_i =V(l_i) \subset \Pn n$, hence points $L_i \in \Pnd n$ in the dual space.
The union of these points is a point in the Hilbert scheme $Hilb_s(\Pnd n)$.  We define the
{\bf v}ariety of {\bf s}ums of {\bf p}owers presenting $F$ as the closure
$$VSP(F,s) = \overline{\{ \{L_1,\ldots,L_s \} \in Hilb_s(\Pnd n) \mid \exists \lambda_i \in \C
: f=\lambda_1 l_1^d+\ldots+\lambda_s l_s^d \} }$$ 
of power sums presenting $f$ in the Hilbert scheme.  Note that taking $d^{th}$ roots of the $\lambda_i$
one could put them into the equations $\l_i$. \par
We study the set of power sum presentations as subsets of natural  Grassmannians.  There may be several 
Grassmannians to choose from, and it is not clear to us whether the compactification in these
Grassmannians always coincide with the Hilbert scheme compactification for a given polynomial.  In
the cases treated here they do. The natural Grassmannians where we embed $VSP(F,s)$ are presented in
section 2.\par
   Both for special and  for general hypersurfaces 
$F=V(f)$ the most basic questions about the schemes $VSP(F,s)$ are:

\bigskip
{\bf Questions:} What is the degree deg$VSP(F,s)$ in case dim$VSP(F,s) = 0$?\par
Is $VSP(F,s)$ irreducible and smooth, in case dim$VSP(F,s) > 0$? 
  \medskip
The schemes $VSP(F,s)$ are contravariants of the hypersurfaces $F \subset \Pn n$ under the action of
$PGL(n+1)$. They attracted an enormous amount of work during the last decades of the $19^{th}$
century, most notably in the work started by [Sylvester 1851], [Rosanes 1873], [Reye 1874] and
[Scorza 1899]. Nevertheless little is known on e.g. the degrees resp. the global structure of these
schemes.  

\bigskip
 Our interest in these schemes arose from  work of Iarrobino [ - 1994] and work of Mukai, 
who proved:
\proclaim 0.6 Theorem {\rm [Mukai 1992]}.  Let $F \subset \Pn 2$ be a general plane quartic. Then
$VSP(F,6)$ is a smooth Fano 3-fold $V_{22}$, i. e. of index 1 and genus 12 with anti-canonical 
embedding of degree 22. Moreover, every $V_{22}$ arises this way.\par

\medskip
 
The fact that plane quartics $F$ are exceptions of Theorem 0.2 lead to interesting geometric 
properties of a $V_{22} \cong VSP(F,6)$: There are 6 conics passing through a general point of
a $V_{22}$, each corresponding to one of the $6$ summands of the power sum. This was observed
classically, e.g. by [Rosanes 1873] and [Scorza 1899]. However the Fano property was not known
classically. Fano overlooked the existence of the $V_{22}$'s in his famous paper [Fano 1937].

{\vskip10pt \bf 0.7 } For the other exceptional cases in Theorem 0.2 there are interesting subvarieties of 
$VSP(F,s)$ by similar arguments. In case of  general quartics $F \subset \Pn 3$ or $\Pn 4$
 or  general cubics $F \subset \Pn 4$ the varieties $ VSP(F,s)$ with $s=10,15$ or $8$ respectively,
are all 5-folds. The main emphasis of this paper lies on the study of $VSP(F,8)$ for cubics in
$\Pn 4$, perhaps the  easiest case among those three. Our main result is this:

\proclaim 0.8 Theorem. Let $F \subset \Pn 4$ be a general cubic.
Then $VSP(F,8)$ is a smooth Fano 5-fold of index 1 and degree 660. More precisely:
Let $\Sp^{10} \subset \Pn {15}$ denote the spinor variety
of isotropic $\Pn 4$'s in the 8-dimensional smooth quadric $Q^8 \subset \Pn 9$. There exists a linear 
subspace $\Pn {10} \subset \Pn {15}$, which depends on $F$, such that $VSP(F,8)$ is isomorphic to 
the variety of lines in the 5-fold $Y=Y(F) := \Pn {10} \cap \Sp^{10} \subset \Pn {15}$.\par

{\vskip10pt} On first sight the occurrence of the spinor variety might look like a surprise. However
it is  less so in view of apolarity and the following result of Mukai:
\bigskip
\proclaim 0.9 Theorem {\rm [Mukai 1995]}.  Let $C \subset \Pn 6$ be a smooth canonical curve of genus
7 and
 Clifford index 3, (e.g. a general canonical curve of genus 7).
Then $C \cong \Pn 6 \cap \Sp^{10}$ for a suitable linear subspace $\Pn 6 \subset \Pn
{15}$.\par   \medskip
Our second result generalizes Mukai's theorem to the case of ``the general empty Gorenstein
subscheme of degree 12 in $\Pn 4$'': Let $R_\Sp$ denote the homogeneous coordinate ring of
$\Sp^{10} \subset \Pn {15}$. Let $A = \C[x_0,\ldots,x_4]/I$ be a graded Artinian Gorenstein ring
with Hilbert function (1,5,5,1). Then we have:

\proclaim 0.10 Theorem.   Let $R_\Sp$ and $A$ be as above. If $A$ is
general, then there exist a regular sequence of linear
forms
$h_0,\ldots,h_{10} \in R_\Sp$
such that $A \cong R_S/(h_0,\ldots,h_{10})$.\par
\medskip

The restriction general here means that the result hold for the general cubic dual polynomial.
The connection to Mukai's theorem is the following: Let $R_C$ be the homogeneous coordinate
ring of a smooth canonical curve $C \subset \Pn 6$. Then for a regular sequence $h_0,h_1 \in
R_C$ of linear forms the Artinian ring $\tilde A =  R_C/(h_0,h_1)$ has Hilbert function
$(1,5,5,1)$. By [Schreyer 1986]  $R_C$ hence $\tilde A$ is syzygy general, iff $C$ has
Clifford index 3.

The connection with Theorem 0.8 is projective duality: The quadric $Q^8 \subset \Pn 9$ has 2 families
of isotropic $\Pn 4$'s, hence two (isomorphic) spinor varieties $\Sp^{10}_{ev} \subset \Pn {15}_{ev}$
and  $\Sp^{10}_{odd} \subset \Pn {15}_{odd}$, where $\Pn {15}_{ev} = \Pnd {15}_{odd}$
naturally. With this identification $\Sp^{10}_{ev}$ is the dual variety $\check \Sp^{10}_{odd}$, cf.
[Ein 1986]. So, if $A^F$ denotes the apolar Artinian Gorenstein ring of a general cubic $F \subset
\Pn 4$, cf. section 1, then $A^F$ has  Hilbert
function (1,5,5,1) and is syzygy general. If $ \Pn {4} \subset \Pn {15}_{ev}$ denotes the linear space
defined by $h_0= \ldots =h_{10} = 0$, where $A^F \cong R_{{\Sp_{ev}}}/(h_0,\ldots,h_{10})$, then 
$\Pn {10} =  {\Pn 4}^\bot \subset \Pn {15}_{odd}$ is the linear space such that 
$Y(F) \cong \Pn {10} \cap \Sp^{10}_{odd}$.

{\vskip10pt \bf 0.11 } {\bf Acknowledgment.} This paper was initiated at a syzygy meeting at Northeastern University 1995,
and completed during the special year of enumerative geometry at Institut Mittag-Leffler 1996/97.   Thanks to Tony Iarrobino and Dan
Laksov and to Stein Arild Str\o mme for helpful discussions on computations.

{\vskip10pt \bf
0.12 } {\bf Notation.} $\C$ denotes an algebraically closed field of characteristic 0. However with
minor modifications all results in this paper hold for arbitrary fields of characteristic zero or
fields of sufficiently large characteristic p, e.g. $char(k) > d$ is necessary even to define the
variety of power sums. All computer experiments in [MACAULAY], which lead us to discover our
results, were done over a finite field.  $\Gr(d,n)$ denotes the Grassmannian of
$d$-dimensional subspaces of
$k^n$, while $\Gr(n,d)$ denotes the Grassmannian of $d$-dimensional quotient spaces of $k^n$.
\medskip We give the numerical information of the minimal free resolution of a graded 
$S=\C[x_0,\ldots,x_r]$-module 
$$ 0 \leftarrow M \leftarrow F_0  \leftarrow F_1 \leftarrow \ldots  \leftarrow F_n  \leftarrow  0$$
with $F_i = \bigoplus_{j \in \Z} \beta_{ij}S(-j)$ in MACAULAY notation, i. e. in the form
$$\matrix{ \beta_{00} & \beta_{11} & \beta_{22} & \ldots & \beta_{n,n} \cr
\beta_{01} & \beta_{12} & \beta_{23} & \ldots & \beta_{n,n+1} \cr
\vdots & \vdots & \vdots & \ldots & \vdots \cr
\beta_{0m} & \beta_{1,m+1} & \beta_{2,m+2} &\ldots &\beta_{n,n+m}. \cr }$$
Note that the integer $m$ is the Castelnuovo-Mumford regularity of $M$. We indicate
$\beta_{ij}$'s which are zero by a $-$. For example the syzygies of the twisted cubic in $\Pn 3$ 
look like: 
$$\matrix{1 & - & - & \cr - & 3 & 2 & ,\cr}$$ 
  there are 3 quadric generators of the ideal
which have 2 linear syzygies.
\vskip20pt

{\bf 1. Apolarity and syzygies}

{\vskip10pt \bf 1.1 }  Consider $S=\C[x_0,\ldots,x_n]$ and 
$T=\C[\partial_0,\ldots,\partial_n]$. $T$ acts on $S$ by differentiation:
$$\partial^{\alpha}(x^{\beta}) =  \alpha!{\beta \choose \alpha } x^{\beta-\alpha}$$
if $\beta \geq \alpha$ and 0 otherwise. Here $\alpha$ and $\beta$ are multi-indices,
${\beta \choose \alpha} = \prod {\beta_i \choose \alpha_i}$ and so on.
In particular we have a perfect pairing between forms of degree d and 
homogeneous differential operators of order d. 
Note that the polar of a form  $f \in S$ in a point $a \in \Pn n$ is given by  $P_a(f)$ for
  $a=(a_0,\ldots,a_r) $ and $P_a  = \sum a_i \partial_i  \in T$.

One can interchange the role of S and T by defining
$$x^{\beta}(\partial^{\alpha}) =  \beta!{\alpha \choose \beta } \partial^{\alpha-\beta}.$$
With this notation we have for forms of degree n
$$P_a^d(f)=f(P_a^d)=d!f(a).$$
Moreover $$f(P_a^m)=0 \iff f(a)=0 \eqno(1.1.1)$$ if $m \geq d$. 
More generally we say that homogeneous forms $f\in S$ and $D\in T$ are {\bf apolar} if $f(D)=D(f)=0$  (According to
[Salmon 1885] the term was coined by Reye). 

\vskip10pt
 Apolarity allows to define Artinian Gorenstein graded quotient rings of $T$ 
via forms: For $f$ a homogeneous form of degree d and $F =V(f) \subset \Pn n$ define
$$ F^{\bot} = f^{\bot} = \{D \in T | D(f)=0 \}$$
and $$A^F = T/F^{\bot}.$$
The socle of $A^F$ is in degree d. Indeed 
$P_a(D(f)) = 0 \hskip3pt \forall P_a \in T_1 \iff D(f) = 0 \hskip3pt or  \hskip3pt D \in T_d$. 
In particular the socle of $A^F$ is 1-dimensional, and $A^F$ is indeed Gorenstein and is called the {\bf apolar Artinian
Gorenstein ring} of $F\subset \Pn n$. 
 
Conversely for a graded Gorenstein ring $A = T/I$ with socle degree d
 multiplication in A induces a linear form $f\colon S_d(T_1) \to \C$ which can be identified with a
 homogeneous
polynomial $f \in S$ of degree d. This proves:

\proclaim 1.2 Lemma {\rm [Macaulay 1916]}. 
The map $F \mapsto A^F$ is a bijection between hypersurfaces $F=V(f) \subset \Pn n$ of 
degree d and graded Artinian Gorenstein quotient rings $A = T/I$ of T with socle degree d.\par 

  Macaulay's result in terms of "inverse systems" is explained more generally in [Eisenbud
1995] Theorem 21.6 and exercise 21.7. The polynomial $f$ is called the {\bf dual socle
generator} or the {\bf dual polynomial} of $A^F$.  Note that the dual polynomial is defined
only up to nonzero scalar.

{\vskip10pt \bf 1.3 } In the following we identify S with homogeneous coordinate ring of $\Pn n$
and T with the homogeneous coordinate ring of the dual space $\Pnd n$. $F=V(f) \subset \Pn n$
denotes always a hypersurface of degree d. We call a subscheme $\Gamma \subset \Pnd n$
$\bf {apolar} $ to $F$, if the homogeneous ideal $I_{\Gamma} \subset F^{\bot} \subset T$.

For example, if $f = l_1^d+\ldots+l_s^d$ and $\Gamma = \{ L_1,\ldots ,L_s\} \subset \Pn n$
the collection of hyperplanes $ L_i =V(l_i) \subset \Pn n$, then $\Gamma$ is apolar to $F =V(f)$.
Indeed $g(l_i^d)$=0 for $g \in I_\Gamma$ by (1.1.1).\par
Conversely, if $\Gamma = \{ L_1,\ldots ,L_s\} \subset \Pnd n$  is apolar to $F$,
then the inclusions $(0)\subset I_{\Gamma}\subset F^{\bot}$ induces inclusions
$$Hom(A_d^F,k)\subset Hom(R_d,k) \subset Hom (T_d,k)$$
where $R=R_{\Gamma}=T/I_{\Gamma}$ is the homogeneous coordinate ring of $\Gamma$.
The linear forms $$\{ D\mapsto D(l_i^d)\}$$
span $Hom(R_d,k)$ so $\{D\mapsto D(f)\}\in Hom(A_d^F,k)$ takes the form
$D(f)=D(\lambda_1l_1^d+\ldots +\lambda_sl_s^d)$ for suitable $\lambda_1,\ldots,\lambda_s\in k$. 
Hence $f= \lambda_1l_1^d+\ldots +\lambda_sl_s^d$  and $\Gamma\in VSP(F,s)$.

Let ${\cal H} \subset Hilb(\Pnd n)$ be a component or a union of components of 
the Hilbert scheme. Then we call the scheme
$$VPS(F,{\cal H}) = \{ \Gamma \in {\cal H} \mid \Gamma \hbox{ is apolar to } F \}$$
the {\bf v}ariety of a{\bf p}olar {\bf s}chemes to $F$ in ${\cal H}$. For $p \in \Q[t]$ we 
abbreviate
$$VPS(F,p) = VPS(F,Hilb_p(\Pnd n)).$$
For $s \in {\bf N} \subset \Q[t]$ we have
$$VSP(F,s) \subset VPS(F,s),$$
This can be a proper inclusion even after taking closures in the Hilbert scheme. 
$VPS(F,p)$ are contravariants of $F$.

{\vskip10pt \bf 1.4 } The consideration of $A^F$ gives a quick explanation, why quartics in
$\Pn 2, \Pn 3 $ and $\Pn 4$ are exceptions of Theorem 0.2: If a quartic is a sum of 5, 9 or 14
powers respectively, then $F^\bot$ contains a quadric. However for a general $F$ the ideal
$F^\bot$ is generated by cubics in all three cases.

In case of cubics $F \subset \Pn 4$ the argument uses syzygies: As T-module an $A^F$ has syzygies
$$ \matrix{ 1 & - & - & - & - & - \cr
                                               - & 10 & 16 & - & - & - \cr
            - & - & - & 16 & 10 & - \cr
              - & - & - & - & - & 1 \cr }$$
for a general $F$. On the other hand for cubics, which are the sum of 7 general powers, $F^\bot$ contains
the ideal $I_C$ of the rational normal curve $C \subset \Pnd 4$ which passes through the 
7 points.
Since the homogeneous coordinate ring of $C$ has syzygies 
$$ \matrix{ 1 & - & - & - & - & -& \cr
                                               - & 6 & 8 & 3 & - & -&, \cr
            }$$
the corresponding Gorenstein ring has syzygies at least
$$ \matrix{ 1 & - & - & - & - & - &\cr
                                               - & 10 & 16 & 3 & - & - &\cr
            - & - & 3 & 16 & 10 & - &\cr
              - & - & - & - & - & 1 &. \cr }$$
So these $F$ are not general.

{\vskip10pt \bf 1.5 } Syzygies of $A^F$ also give a quick uniform proof of the assertions in 0.4:
For a binary form $f$ of odd degree $d=2k-1$ the apolar Artinian Gorenstein ring is a complete
intersection $$A^F \cong \C[\partial_0,\partial_1]/(a,b) $$
with ${\rm deg} a +{\rm deg} b =2k+1$, say ${\rm deg} a < {\rm deg} b$. For a general $f$ we have ${\rm deg} a = k$ and the 
k roots give the unique set $l_1,\ldots,l_k \in \C[x_0,x_1]$ such that 
$f = \lambda_1l_1^d+\ldots+\lambda_kl_k^d$.
 Similarly, of course, for a general form of even degree $d=2k$, the apolar Artinian Gorenstein
 ring is a complete intersection
$$A^F \cong \C[\partial_0,\partial_1]/(a,b) $$
with ${\rm deg} a={\rm deg} b =k+1$. Thus $VSP(F,k+1)\cong \Pn 1$. 

In Hilbert's case $A^F$ is Gorenstein of codimension 3 and the structure theorem of 
Buchsbaum-Eisenbud applies. In particular the number of generators of $F^\bot$ is odd.
 For a general quintic $A^F$ has syzygies
$$ \matrix{ 1 & - & - & - &\cr
                - & - & - &- & \cr
                - & 4 & 1 &- &\cr
            - & 1 & 4 & - &\cr
                - & - & - &- & \cr
              - & - & - &1 &, \cr }$$
with the middle matrix $\phi$ is skew symmetric, its Pfaffians generate the ideal. The 4 linear entries of the first column are dependent,
since $\C[\partial_0,\partial_1,\partial_2] $ has only 3 linearly independent linear forms.
Thus after column and row operations we may assume
$$ \phi = \pmatrix{ 0 & 0 & \partial_0 & \partial_1 & \partial_2 & \cr
0 & 0 & q_{23} & q_{24} & q_{25} & \cr
-\partial_0 & -q_{23} & 0 & q_{34} & q_{35} & \cr
-\partial_1 & -q_{24} & -q_{34} & 0 & q_{45} & \cr
-\partial_2 & -q_{25} & -q_{35} & -q_{45} & 0 & \cr} ,$$
for a general $F$. The $2 \times 2$ minors of the block
$$\pmatrix{\partial_0 & \partial_1 & \partial_2 \cr
q_{23} & q_{24} & q_{25}  \cr}$$
are among the Pfaffians of $\phi$ and generate the ideal of the unique 7 points in $\Pnd 2$.

Finally for the pentahedral theorem we note that $A^F$ for a general cubic $F \subset \Pn 3$ has
syzygies
$$ \matrix{ 1 & - & - & - &  - &\cr
                                               - & 6 & 5 & - &  - &\cr
            - & - & 5 & 6 & - &\cr
              - & - & - & - & 1 & \cr } \cong \quad 
 \matrix{ 1 & - & - & - & \cr
                                               - & 5 & 5 & - &\cr
           - & - & - & 1 & \cr } \otimes \quad \matrix{ 1 & - & \cr - & 1 &, \cr}$$
cf. [Kustin, Miller 1985] for the decomposition.
The 5 Pfaffians among the 6 quadrics define the ideal of 5 distinguished points. (This decomposition of syzygies is well-known in case
of canonical curves $C \subset \Pn 5$ of genus 6: A general genus 6 canonical curve $C$ is a complete
intersection of an unique Del Pezzo surface of degree 5 and a non unique hyperquadric, cf. [ACGH
1985], [Schreyer 1986], [Mukai 1988].)

{\vskip10pt \bf 1.6 } In the cases above one can give examples of forms which need more summands
than the general form for a power sum presentation easily: Just take the unique finite length subscheme
non-reduced.
 
For binary forms the syzygy approach solves the question posed in 0.3 and 0.5. E.g.
every binary form of degree d is sum of $\tilde s = d$ or less  $d^{th}$ powers. An example of a form
which achieves this bound is $f=xy^{d-1}$. Here 
$$A^F = \C[\partial_x,\partial_y]/(\partial_x^2,\partial_y^d).$$
Since every element of $F^\bot$ of degree $\le d-1$ has a double root $\tilde s = d.$ For
$\partial_x^d-\partial_y^d \in F^\bot$ the corresponding power sum is
$$d^2xy^{d-1} = \sum_{j=1}^d ({\zeta_{d^2}^jx+\zeta_{d^2}^{jd+j} y} )^d $$ 
where $\zeta_{d^2}$ is a primitive $d^2$-th root of unity.

{\vskip10pt}  Mukai's result $VSP(F,6) \cong V_{22}$ for a general plane quartic $F \subset
\Pn 2$, is obtained with syzygies.  The approach is similar for general planes curves of even degree
up to 8. 
\proclaim 1.7 Theorem. Let $F$ be a general plane curve of degree 
$d=2n-2$, $2\leq n\leq 5$,
then  $$VSP(F,{{n+1}\choose 2}) \cong \Gr  (n,V,\eta) = \{ E \in \Gr  (n,V) \mid \Lambda^2 E
\subset ker(\eta) \}$$ where $V$ is a $2n+1$-dimensional vector space and $\eta$ is a net of
alternating forms $\eta \colon \Lambda^2 V\to \C^3$ on $V$.\smallskip
i) {\rm [Mukai 1992]} When $F$ is a smooth plane conic section, then $VSP(F,3)$ is a Fano 3-fold of index 2 and degree 5 in $\Pn 6$.\smallskip
ii) {\rm [Mukai 1992]} When $F$ is a general plane quartic curve, then $VSP(F,6)\cong V_{22}$ (cf.
{\rm 0.6}).\smallskip iii) {\rm [Mukai 1992]} When $F$ is a general plane sextic curve, then
$VSP(F,10)$ is isomorphic to polarized $K3$-surface of genus 20.\smallskip iv) When $F$ is a
general plane octic curve, then $VSP(F,15)$ is finite of degree 16 , i. e. consists of 16
points.\par $Proof.$  A ternary form $f$ of even degree $d=2n-2$ is called 
non-degenerate, if $F^\bot$ contains no elements of degree $n-1$. 
 Here, as always, $F=V(f)$ denotes the corresponding plane curve
of degree $d$. For non-degenerated $f$ the Artinian Gorenstein ring $A^F$ has syzygies

$$ \matrix{ 1 & 
- & - &  - &\cr
           - &  - & - &  - &\cr
\vdots & \vdots &\vdots & \vdots&\cr
          - &  2n+1 & 2n+1 & - &\cr
 \vdots &\vdots &\vdots &\vdots&\cr                                                                  
                - &  - & - &  - &\cr
              - &  - & - & 1 & ,\cr }$$
because by the symmetry of the Hilbert function of $A^F$ there are
$2n+1 = { n+2 \choose 2} - {n \choose 2}$ generators of degree $2n$ in $F^\bot$
and by the symmetry of the resolution there are no generators in higher degree.
Conversely a general skew symmetric $2n+1 \times 2n+1$ matrix $\phi$
 of linear forms
defines via its Pfaffians an Artinian ring with socle in degree $d=2n-2$, 
whose form is non-degenerate.

Interpreting $\phi$
 as a net of alternating forms $\eta \colon \Lambda^2 V\to \C^3$
 on a $2n+1$-dimensional vector space $V = (F^\bot)_n^*$ one 
 obtains:
 \medskip

\proclaim 1.8 Lemma.  Let f be a non-degenerate ternary form of degree
$d=2n-2, n=2,3,4,5$ and $\eta \colon V \to \bC^3$ the corresponding
alternating net. Then 
$$VPS(F,{n+1 \choose 2}) = \bG(n,V,\eta)$$\par

$Proof.$ Let $X \in \bG(n,V,\eta)$ be a point, and let
  
$$\phi = \pmatrix{ B & A \cr
-A^t & 0 \cr}$$
be the corresponding decomposition of the syzygy matrix with a
$n+1 \times n+1$ skew-symmetric matrix $B$,  a $n+1 \times n$ matrix $A$ and
a $n \times n $ block of zeros. Among the Pfaffians of the
syzygy matrix are the $n \times n$ minors of $A$. We claim that the
minors define a scheme of length $n+1 \choose 2$ in $\bP^2$. By
Hilbert-Burch [Eisenbud, Thm 20.15 ], this is the case unless the $n+1$
minors $f_0,\ldots,f_n$ have a common factor $h$ of degree $k > 0$. 
We will derive a contradiction by studying the syzygies of the ideal
$J=(g_0,\ldots,g_n)$ with $f_i = hg_i$. $g_0,\ldots,g_n$ are
linearly independent, since $f_0,\ldots,f_n$ have this properties.
Their free resolution has shape:
$$0 \leftarrow T/J \leftarrow T \leftarrow (n+1)T(-n+k) \leftarrow
nT(-n+k-1)\oplus T(-2n+2k-1) \leftarrow T(-2n+k-1) \leftarrow 0.$$
Indeed syzygies among the $g$'s are syzygies among the $f's$, hence
linear combinations of the columns of $\phi$ with  n zero
components. Since $ker -A^t = (g_0,\ldots,g_n)^t$ the syzygies are
generated by the columns of
$A$ and the vector $B(g_0,\ldots,g_n)^t$. The second syzygy
is the vector $(f_{n+1},\ldots,f_{2n},h)^t$. Here $f_{n+1},\ldots$
denote the remaining Pfaffians.
It follows
$$ {\rm deg}T/J = {n+1 \choose 2} - (2n+1)k +k^2$$
which gives a negative number in the possible range of $(n,k)$,
unless $k=1$ and $n=3,4,5$. In these remaining cases we consider
the resolution of the last matrix transposed. The n forms of degree n
modulo the linear form h, generate an ideal of codimension 2 in T/h,
hence have a determinantal resolution
$$ 0 \leftarrow T/J_1 \leftarrow T/h
\leftarrow n T/h(-n) \leftarrow (n-2)T/h(-n-1) \oplus T/h(-n-2)
\leftarrow 0,$$
where $J_1 = (h,f_{n+1},\ldots,f_{2n}),$
and over $T$ the resolution is
$$0 \leftarrow T/J_1 \leftarrow T \leftarrow
T(-1) \oplus nT(-n) \leftarrow (2n-2)T(-n-1)\oplus T(-n-2) \leftarrow
(n-2)T \oplus T(-n-3) \leftarrow 0.$$
Consider the sheafivication of this complex and the middle syzygy sheaf
${\cal F}$.
$(g_0,\ldots,g_n)$ generate the kernel of
$$  {\cal F}{\buildrel \alpha \over \leftarrow} (n+1){\cal O}(-n-1)  \leftarrow
{\cal O}(-2n) \leftarrow 0.$$
However from the presentation of ${\cal F}$
$$0 \leftarrow {\cal F} \leftarrow (2n-2){\cal O}(-n-1) \oplus {\cal
O}(-n-2)
\leftarrow (n-2){\cal O}(-n-2) \oplus {\cal O}(-n-3) \leftarrow 0$$
we obtain a short exact sequence
$$0 \leftarrow {\cal O}_H(-n-2) \leftarrow {\cal F} \leftarrow {\cal G}
\leftarrow 0.$$
So $Im \alpha \subset {\cal G}$ and 
$c_1(ker \alpha) \ge -(n+1)^2-c_1({\cal G})=-(n+1)^2-c_1({\cal
F})+1=-2n+1$
a contradiction to the degrees of the $g's$.

Thus 
$$VPS(F,{n+1 \choose 2}) \supseteq \bG(n,V,\eta).$$  For the other inclusion,
 consider a finite subscheme $\Gamma\in VPS(F,{n+1 \choose 2})$.
Since $F$ is nondegenerate,  $\Gamma$ imposes independent conditions on forms of degree $n-1$.
By Hilbert-Burch $\Gamma$ has syzygies:
$$ 0 \longleftarrow R_{\Gamma} \longleftarrow T\longleftarrow (n+1)T(-n)
\buildrel \psi \over 
\longleftarrow nT(-n-1) \longleftarrow 0. \eqno(1.8.1)$$
By apolarity the ideal $J_{\Gamma} \subset F^{\bot}$. Hence we have a
sequence
$$ 0 \longleftarrow A^F  \longleftarrow R_{\Gamma} \longleftarrow F^{\bot}/J_{\Gamma}
\longleftarrow 0. \eqno(1.8.2)$$
Since $A^F$ and $R=R_{\Gamma}$ have Hilbert functions 
$$(1,3,\ldots {n+1 \choose 2}, {n \choose 2} \ldots,3,1)$$ and
$$(1,3,\ldots ,{n+1 \choose 2}, {n+1 \choose 2}, {n+1 \choose 2},\ldots)$$
respectively, $F^{\bot}/J_{\Gamma}$ has $n$ generators of degree $n$ with $n+1$ linear relations:
$$0 \longleftarrow F^{\bot}/J_{\Gamma} \longleftarrow nT(-n) \longleftarrow (n+1)T(-n-1).$$
The minors of the presentation matrix are contained in the annihilator,
which is $J_{\Gamma}$.
Hence, this matrix is $\psi^t$ again, and $F^{\bot}/J_{\Gamma} \cong \omega_R(-2n+2)$. A
mapping cone between
the complex (1.8.1) and its dual over the sequence (1.8.2), gives syzygies
of $A^F$ with the
desired block structure.$\slutt$\medskip 

{\it End of the proof of} (1.7).
General curves $F$ correspond via dual socle generators to general
alternating matrices. Hence for a general $F$ the variety
$\bG(n,V,\eta)$ is smooth of expected dimension as zero loci of a
general section of a 
globally generated bundle.  In fact, let $Q$ be the universal rank $n$ quotient bundle on
 $\bG=\bG(n,V)$ with chern classes $c_i=c_i(Q)$.
Then $X=\bG(n,V,\eta)$ is the zero locus of a section of $3\wedge^2Q$, so 
$${\rm deg}X=c_1^m\cdot c_{3{n\choose 2}}(3\wedge^2Q),$$ where 
$$\eqalign {m=&{\rm 
dim}X=n(n+1)-3{n\choose 2}={1\over 2}n(5-n)\cr =&3,3,2,0\quad {\rm 
for}\quad n=2,3,4,5.\cr} $$  A  calculation in the 
Chow ring of the respective $\bG$'s give the degrees $5,22,38$ and 
$16$, respectively, for $X$.

 Consider now the incidence correspondence

$$\{ (\Gamma,F) \in Hilb^o({n+1 \choose 2},\bP^2) \times \mid (2n-2)L
\mid^o \mid \Gamma \hbox{ is apolar to } F \}$$
where $^o$ indicates that we consider only the open parts of
$\Gamma$'s, respectively $F$'s, which impose independent conditions
on forms of degree n-1. The incidence correspondence is irreducible,
since its fibers over $Hilb^o$ are irreducible of constant dimension. 
The fibers of the other projection are the $\bG(n,V,\eta)$'s by the
lemma. Since smooth $\Gamma$ are dense in the Hilbert scheme, smooth $\Gamma$ are
dense in every component of $\bG(n,V,\eta)$ for a general $F$. 
This proves (iv). 

In the other cases
all components have the same degree and have induced canonical
bundle by the irreducibility of the correspondence. 

The tangent bundle $T_{\bG}$ has first chern class 
$c_1(T_{\bG})=(2n+1)c_1$, while the normal bundle of $X$ in $\bG$, the 
bundle $N_X=3\wedge^2Q$, has first chern class $3(n-1)c_1\cap X$.  
Therefore the chern class of the canonical bundle on $X$ is 
$(n-4)c_1\cap X$. 
Hence the components are Fano $3$-folds of index $2$
in case deg$F = 2$, Fano $3$-folds of index $1$ in case deg$F = 4$, and 
$K3$- or abelian 
surfaces in case deg$F = 6$. The sums of the degrees of the components
are $5, 22=2\cdot 11$ and $38=2\cdot 19$ respectively. Since $5,11$ 
and $19$ are prime numbers
the fiber is irreducible for a general $F$.
This proves (i) and (ii).

In case (iii) it remains to show that $X$ is $K3$. The total chern 
classes of $T_{\bG}$ and $N_X$ are
$$c(T_{\bG})=1+9c_1+39c_1^2+c_2+\ldots\quad {\rm and}\quad  c(N_X)=(1+9c_1+36c_1^2+6c_2)\cap X,$$
respectively.  Thus $T_X$ has total chern class
$c(T_X)=c(T_{\bG})/c(N_X)=1+3c_1^2-5c_2$.  Since the degree of $X$ is $c_1^2\cap X=38$, 
we get $c_2(T_X)=114-5c_2\cap X$,
 i.e. $X$ can only be K3 (with $c_2\cap X=18$).
$\slutt$\medskip 

 \medskip This proves the easy part of Mukai's theorem for plane quartics. For the difficult
part, that every $V_{22}$ arises this way, cf. [Mukai 1992]. \par

{\vskip10pt \bf 1.9 } To describe $VSP(F,8)$ for a general cubic $F \subset \Pn 4$ we look for
subcomplexes 
$$ \matrix{ 1 & - & - & - & - & \cr
              - & 7 & 8 & - & - & \cr
            - & - & 3 & 8 & 3 & , \cr
             }$$
corresponding to 8 points $\Gamma \subset \Pnd 4$, of the syzygies of $A^F$
$$ \matrix{ 1 & - & - & - & - & - &\cr
              - & 10 & 16 & - & - & - &\cr
            - & - & - & 16 & 10 & - &\cr
              - & - & - & - & - & 1 &.\cr }$$
Thus we have embeddings $VSP(F,8) \hookrightarrow \Gr  (7,10)$ or $VSP(F,8) \hookrightarrow \Gr (8,16)$.
However unlike Mukai's case, we could not find a description of the compactifications of those subschemes, which 
would give  e.g. a correct dimension estimate. Our method is more involved.
\vskip20pt

{\bf 2. Duality, projections and secants}

{\vskip10pt \bf 2.1 } The syzygy approach of the previous section has a geometric counterpart in terms
of duality and projections which will clarify our solution of the case of a cubic
in $\Pn 4$.  In the next section we shall see how this more geometric approach and the syzygy
approach both apply to the cases of plane curves of degree 7 and quadric
surfaces.\par 
With notation as above let $f\in S_d$ be a homogeneous
 form of degree $d$, and assume $f = l_1^d+ \ldots + l_s^d$ for some $s$ and 
some linear forms $l_i$ i. e. $\Gamma = \{ L_1,\ldots ,L_s\} \subset \Pnd n$
is apolar to $F=V(f)$.
  The first simple but crucial observation is that for any homogeneous $D\in T$ of degree
$0\leq d^{\prime}\leq d$, we get  $$D(f)=\sum_{i=1}^s\lambda_i^Dl_i^{d-d^{\prime}},$$
with $\lambda_i^D\in \C$. 
Thus $\Gamma$ is apolar to $F_D=V(D(f))$, i. e. to all multiple partials of $f$.

{\vskip10pt \bf 2.2 }
  Now fix an $e\leq d$ and let
$$F_e^{\bot}=\{D\in T_e|D(f)=0\}\subset F^{\bot}.$$
Then  $$F_e^{\bot}=\{D\in T_e|D(D^{\prime}(f))=0 \quad
{\rm for\quad all}\quad D^{\prime}\in T_{d-e}\}\subset T_e,$$
so $F_e^{\bot}$ is the dual space to the space
$$F_e=\{D^{\prime}(f)|D^{\prime}\in T_{d-e}\}\subset S_e$$ of multiple partials of $f$.
 We consider the map 
$$ \pi_e^F:\Pnd n\to {\bf P_{\ast}}(F_e^{\bot})=\Pn {n_e}$$
 defined by $F_e^{\bot}$, with $n_e={\rm dim}F_e^{\bot }-1$. By construction $\pi_e^F$ is
 the composition of the $e$-uple embedding of $\Pnd n$ and the {\bf projection from} 
the {\bf partials}
 of $f$ of order $d-e$ considered as points in $\Pnd {{{n+d}\choose d}-1}$.
{\vskip10pt \bf 2.3}  If $\Gamma = \{ L_1,\ldots ,L_s\} \subset \Pnd n$
 is apolar to $F_e$, then $I_{\Gamma}(e)\subset F_e^{\bot}$ and the span of $\pi_e^F(\Gamma)$
has dimension  $$d_{\Gamma}(e)={\rm dim}F_e^{\bot}-{\rm dim}I_{\Gamma}(e)-1.$$ 

Therefore, when $F_e^{\bot}$ has no basepoints, $\pi_e^F(\Gamma)$  is contained in a fiber
 of $\pi_e^F$ over a linear space of dimension $d_{\Gamma}(e)$.  When  $\pi_e^F$ is an
 embedding the apolar set $\Gamma$ determine a $s-$secant  $\Pn {d_{\Gamma}(e)}$ to 
$\pi_e^F(\Pnd n)$ this way.
 
{\vskip10pt \bf 2.4}  If $d_{\Gamma}(e)=d_e(s)$ is independent of $\Gamma\in
VSP(F,s)$, then $$\Gamma\mapsto <\pi_e^F(\Gamma)>$$ defines a map
 $$\rho:VSP(F,s)\to \Gr (d_e(s)+1,n_e+1).$$
This occurs naturally when the $(d-e)$-uple partials of $f$ are linearly
independent and the points of $\Gamma$
impose independent conditions on forms of degree $e$ for any $\Gamma\in
VSP(F,s)$. In this case $n_e={{n+d}\choose {d}}-{{n+d-e}\choose {d-e}}$ and  
$$d_{\Gamma}(e)=d_e(s)=s-{{n+d-e}\choose {d-e}}-1.$$ 

 In general $\rho$ is a rational map, when it is
birational it gives a birational model of $VSP(F,s)$.  Of course, this model may not coincide with the Hilbert
scheme compactification of (0.5), and in general we see no reason why it should.\par

In the cases treated here we choose the minimal $e$ such that $d_{e}(s)\geq 0$.  This number
coincides with the highest degree of the generators of $F^{\bot}$ and the homogeneous ideal
$I_{\Gamma}$ of $\Gamma\in VSP(F,s)$. Since $A^F$ is Artinian, $F_e^{\bot}$ cannot have base points, and the map 
 $\rho:VSP(F,s)\to\Gr (d_e(s)+1,n_e+1)$ is a morphism.  In each case treated here we check in fact that it is an
embedding of $VSP(F,s)$.\par
 By construction the image of $\rho$ is a subvariety of the variety of $s$-secant spaces to
$\pi_e^F(\Pnd n)$ of dimension $d_e(s)$, and in general it is proper.\par 
 
In fact, the image $\pi_e^F(L)$ of any line $L\subset \Pnd n$ have a span of
 dimension at most $e$.  So as soon as $e\leq d_e(s)$ any scheme of length $s$ on the line will contribute
to the variety of $s$-secant spaces to
$\pi_e^F(\Pnd n)$ of dimension $d_e(s)$, while such schemes are not apolar to $F$.\par 
The general criterion for a subscheme $Z$ to be apolar to $F$ can be weakened to a useful sufficient
criterium: Consider a scheme $Z\subset \Pnd n$ of length $s$ such that the span of
$\pi_e^F(Z)$ has dimension $d_e(s)$.  Let $I_Z$ be the homogeneous ideal of $Z$.  If $Z$ impose $s$ conditions
on forms of degree $e$, then $I_Z(e)\subset F_e^{\bot}$.  If furthermore  $I_Z$ is generated by forms of
degree $e$, then $I_Z\subset F^{\bot}$ and $Z$ is apolar to $F$. \par

{\vskip10pt \bf 2.5}  The projection from partials is well illustrated by the case of the cubic
consisting of three lines in the plane. Let $f=x_0x_1x_2, \quad F=V(f)$.  Then 
$$F_2^{\bot}=\{D\in T_e|D(f)=0\}=<\partial_0^2,\partial_1^2,\partial_2^2>,$$
and $$\pi_2^F:\Pnd 2\to\Pn 2$$ is a $4:1$-map.
If $\Gamma\in VSP(F,4)$, then $I_{\Gamma}$ is generated by two quadrics, so $\Gamma$ is a fiber of $\pi_2^F$.
Conversely any general fiber $Z$ of $\pi_2^F$ is a collection of 4 distinct points, 
the complete intersection of two quadrics.
Therefore $Z$ is apolar to $F$.  We conclude that $VSP(F,4)\cong \Pn 2$.\par
Note that this approach works for the general plane cubic so $VSP(F,4)\cong \Pn 2$ 
for any general cubic $F\subset \Pn 2$. \par
{\bf Remark.} We get the Hilbert scheme compactification in this case since the map 
$\pi_2^F$ is finite and every fiber have length 4.

{\vskip20pt}
{\bf 3. Another example}
\bigskip
In this section we illustrate the syzygy approach and the duality and projection approach 
applied to the case of plane curves of degree 7.  
\proclaim 3.1 Theorem. {\rm [Dixon, Stuart 1908]}.  Let $F\subset \Pn 2$ be a
general  curve of degree 7.  Then
$VSP(F,12)$ consists of 5 distinct points, i. e. the form of degree 7 defining $F$ 
have precisely 5 distinct presentations
as the sum of 12 powers of linear forms.\bigskip

We consider first the apolar graded Artinian Gorenstein ring
$A^F=\C[\partial_0,\partial_1,\partial_2]/F^{\bot}$. It is a certain "Veronese" quotient
of the homogeneous coordinate ring of a smooth Del Pezzo surface $D$ of degree 5 in $\Pn 5$. 
More precisely let $I_D\subset\C[y_0,\dots,y_5]$ be the homogeneous ideal of a smooth Del Pezzo
surface in $\Pn 5$ :
\proclaim 3.2 Proposition.  Let $A^F$ be as above. Then there is a Veronese embedding of $\Pnd 2\to
\Pn 5$ with corresponding ring homomorphism  $\varphi: \C[y_0,\dots,y_5]\to
\C[\partial_0,\partial_1,\partial_2]$ such that 
$A^F=\C[\partial_0,\partial_1,\partial_2]/\varphi(I_D)$ where $\varphi(I_D)$  denotes the ideal
generated by the image of $I_D$ under $\varphi$.\par
$Proof$.  For a general $F$ the Artinian Gorenstein ring $A^F$ has the maximal
possible Hilbert function (1,3,6,10,10,6,3,1). So there are 5 quartic 
generators in $F^\bot$. The syzygies are

$$ \matrix{ 1 & - & - & - &\cr
- & - & - & - & \cr
- & - & - & - & \cr
- & 5 & -& - &\cr
- & - & 5 & - &\cr
- & - & - & - & \cr 
- & - & - & - & \cr
- & - & - & 1 &. \cr }$$
because those without relations in degree 6 are the most general. 
The syzygy matrix has 10 off diagonal quadratic entries, and it is an open 
condition on $F$ that these generate all 6  linearly independent quadrics in
$\Pn 2$. Thus we obtain a 
Veronese embedding $\Pn 2 \hookrightarrow \Pn 5 \subset \Pn 9$
in which we have the Del Pezzo surface $\Pn 5 \cap \Gr (2,5)$. 
The two surfaces do not intersect because $A^F$ is Artinian. 
Conversely given a Veronese embedding 
$\varphi \colon \C [y_0,\ldots,y_5] \to \C [\partial_0,\partial_1,\partial_2]$
and a Del Pezzo surface $D$ such that $\Pn 2 \cap D = \empty$,
then $\phi(I_D)$ defines an Artinian Gorenstein ring with desired syzygies
$\slutt$\medskip

{\it Proof of} (3.1): Recall from the proof of the proposition that $F^{\bot}$ is generated 
by 5 quartics, they
are restrictions of the quadrics defining a Del Pezzo surface $D$ to the Veronese embedding 
$V$ of $\Pnd 2$.
The map $\pi_4^F:\Pnd 2\to \Pn 4$ defined by these 5 quartics is therefore a morphism.  
A set of 12 points $\Gamma\subset\Pnd 2$
is apolar to $F$ when they are collinear under this map, and the ideal $I_{\Gamma}$ is 
generated by three quartics, cf. 2.4. 
Now, the Del Pezzo is a linear section of the Grassmannian $\Gr (2,5)\subset \Pn 9$.  
The 5 quadrics generating the ideal of $\Gr (2,5)$
define a rational map $\pi:\Pn 9 \to\Pn 4$.  The preimage of a line is a twisted cubic 
6-fold scroll with vertex a plane, it lies in a
tangent hyperplane to $\Gr (2,5)$.  Similarly the preimage of a line for the map defined by
the quadrics in $I_D$ are restrictions of
these twisted cubic scrolls to $\Pn 5$.  The general one meets $\Pn 5$ in a twisted 
cubic surface, which does
not meet $V$.
The tangent hyperplanes to $\Gr (2,5)$ are naturally parameterized by $\Gr (5,3)$ in
the dual space $\Pnd 9$.  Dual to $\Pn 5\subset \Pn 9$
is a $\Pn 3\subset \Pnd 9$.  This $\Pn 3$ meets $\Gr (5,3)$ in 5 points.  These 5 points 
correspond to the 5 twisted cubic threefold scrolls
that contain $D$.  They each intersect $V$ in a scheme of length 12 which become collinear 
under the map to
$\Pn 4$. The ideal of these
points is generated by the three minors of a $2\times 3$ matrix with quadratic entries, so 
these 5 schemes of length 12 are precisely
the elements of $VPS(F,12)$.  When the Veronese embedding is general with respect to the 
Del Pezzo surface $D$,
these apolar schemes
are all smooth such that $VPS(F,12)=VSP(F,12)$, and the theorem follows. $\slutt$
\medskip

Together with (1.5),(1.7) and (2.5) this result treats general plane curves of degree $\leq 8$.  
For higher degrees the problems seems more difficult.

\vskip20pt
{\bf 4. Cubic threefolds, first properties}

{\vskip10pt \bf 4.0 }  Let $f\in S=\C[x_0,\dots,x_4]$ be a general cubic
polynomial, and $F=V(f)\subset \Pn 4$, and let $T=\C[\partial_0,\dots,\partial_4]$.  In the
remaining sections we investigate $VSP(F,8)$.  Chronologically we started to study its local
geometry.  We show in this section that there are 8 conic sections through each point in
$VSP(F,8)$.  In the next section we study the syzygies of 8 general points in $\Pn 4$, and
in particular describe the pencil of elliptic quintic curves through 8 general points.  
For a general cubic, we then started to look for the variety of apolar elliptic quintic
curves.  With MACAULAY we computed examples of apolar Artinian Gorenstein rings for cubics
$f$ with general syzygies, and found a rank 10 quadratic relation among the apolar
quadrics to $f$.  The quadrics defining an apolar elliptic quintic curve through an apolar subscheme
of length 8, form an isotropic $\Pn 4$ for the corresponding quadratic form. Thus we look for a
certain subvariety $Y$ of the set of isotropic $\Pn 4$'s of a smooth quadric $Q$ in $\Pn 9$
parametrizing apolar elliptic quintic curves. 
 The isotropic
$\Pn 4$'s are parameterized by a 10-dimensional spinor variety $\Sp$, which
plays a crucial role in our investigation.  We introduce $\Sp$ in section 6 with its
natural half spin embedding in $\Pn {15}$, and show that the apolar Artinian Gorenstein ring
of a general cubic $F$ is an Artinian quotient of the homogeneous coordinate ring of this
spinor variety in section 7. In section 8 we then use the facts that the spinor variety is
isomorphic to its dual, and that the quadrics generating its homogeneous ideal define a map
of $\Pn {15}$ onto the smooth quadric $Q$ in $\Pn 9$, to find the variety of
apolar quintic elliptic curves, and prove the main result on $VSP(F,8)$.  In the last
section we show that $VSP(F,8)$ is a Fano variety of
index 1 and degree 660.

\vskip10pt {\bf 4.1}
First we construct a rational map of $VSP(F,8)$ into the Grassmannian of $\Gr (3,10)$ of
planes in $\Pn 9$. The apolar quadrics $$F_2^{\bot}=\{D\in T_2|D(f)=0\},$$ define in the
notation of section 2, a map
$$\pi_2^F:\Pnd 4\to \Pn 9(={\bf P_{\ast}}(F_2^{\bot})),$$
the composition of the Veronese embedding $\Pnd 4\to\Pnd {14}$ with the projection from the
partials of $f$ considered as points in $\Pnd {14}$. For general $F$ this map is an embedding as we
shall see later (cf. 8.3). We denote the image by $V_F$.
\par
As in (2.4), if $\{L_1,\ldots ,L_8\}\in VSP(F,8)$ then 
$$<\pi_2^F (L_1,\ldots ,L_8)>={\rm span}\{\pi_2^F(L_1),\dots ,\pi_2^F(L_8)\}$$ is only a plane
in $\Pn 9$.  It is an 8-secant plane to $\pi_2^F(\Pnd 4)=V_F\subset \Pn 9$.
Thus we have defined a map of $VSP(F,8)$ into $\Gr (3,10)$ as promised.  In section 8 we
show that this map is in fact an embedding. For the purpose of this section we identify
$VSP(F,8)$ with its image in $\Gr (3,10)$.\par 
Notice that all 8-secant planes to $V_F$ need not come from points on $VSP(F,8)$. In fact
any line in $\Pnd 4$ is embedded as a conic on $V_F$, but 8 points on it does not belong to
$VSP(F,8)$.   

\vskip10pt {\bf 4.2}  Recall from $(1.4)$ that for cubics
$F=V(f)\subset\Pn 4$ where  $f$ is a sum of 7 general powers, called a
degenerate cubic, the unique rational normal quartic curve $C\subset \Pnd 4$
passing through the 7 points is apolar to $F$.  By the 3-uple embedding of $\Pnd 4$ into
$\Pn {34}$, the curve $C$ has degree 12 and spans a $\Pn {12}$ which contain the
cubic $f$.  This $\Pn {12}$ can be interpreted as the space of binary forms of degree 12,
and the fact, from $(1.5)$, that $VSP(f_{bin}, 7)\cong\Pn 1$ when one considers $f$ as such a
binary form $f_{bin}$, means that $C$ has a $\Pn 1$ of $7$-secant $\Pn 6$'s passing through
$f$. In the interpretation of $f$ as a cubic the conclusion is that $VSP(F,7)$ contains a
$\Pn 1$.  
On the other hand the quadrics
through $C$ are apolar to $F$. The apolar quadrics
 $$F_2^{\bot}=\{D\in
T_2|D(f)=0\},$$ define a map $$\pi_2^F:\Pnd 4\to
\Pn 9(={\bf P_{\ast}}(F_2^{\bot})).$$   Recall that a rational normal quartic curve is defined by the
$2\times 2$ minors of a $2\times 4$ matrix with linear entries.  Among these minors there is
the Pl\"ucker quadric relation, a rank 6 quadric. Therefore $\pi_2^F(\Pnd 4)\subset
\Pn 9$ is contained in a rank 6 quadric, and $C_F=\pi_2^F(C)$ spans the
vertex $\Pn 3$ of this quadric. Any 7 points apolar to $F$ become collinear on
$C_F$.  Thus $C_F$ is a curve of degree 8 in $\Pn 3$ with at least a $\Pn 1$ of $7$-secant
lines.  This, of course, occurs precisely when $C_F$ is a curve of type $(1,7)$ on a smooth
quadric surface.  Thus for general degenerate cubic $F$, the variety
$VSP({F_{\lambda}},7)\cong\Pn 1$, identified by one of the rulings of a quadric in
a $\Pn 3$.\par
{\vskip10pt \proclaim 4.3 Proposition. There are 8 conic sections through the general point
on $VSP(F,8)\subset \Gr (3,10)$ corresponding to the 8 summands of the power sum presentation.\par
$Proof$.  Let $\{L_1,\ldots ,L_8\}\in VSP(F,8)$ be a general point with $L_i=V(l_i)$. Then
$f=l_1^3+\ldots +l_8^3$.  Thus $f_i=f-l_i^3$ is a degenerate cubic, and the union of any
element of $VSP(F_i,7)$ with $L_i$ belongs to $VSP(F,8)$.  From (4.2) we get altogether 8 
$\Pn 1$'s through $\{L_1,\ldots ,L_8\}$.  Now to show that these $\Pn 1$'s are conic sections in the
image of $VSP(F,8)$ in $\Gr (3,10)$ we consider again the rational normal curve apolar to $F_i$.
The image of this curve on $V_F$ has degree 8 and has a pencil of 7-secant planes all
passing through the point $L_i$.  The projection from $L_i$ therefore have degree 8 and
have a pencil of 7-secant lines, which must form a quadric surface like above.  Therefore
the 7-secant planes form a pencil of planes of a rank 4 quadric in a $\Pn 4$. But such a
pencil clearly form a conic section in $\Gr (3,10)$, so the proposition follows.$\slutt$ 

\vskip20pt
{\bf 5. Syzygies of 8 general points in $\Pn 4$}\bigskip

 In this section we describe the family of elliptic quintic normal curves through 8 general
points in $\Pn 4$.  
{\vskip10pt \bf 5.1 }Let $\Gamma = \{p_1,\dots,p_8\} = D_1 \cap D_2 \subset S
\subset \Pn 4$ be the 8  intersection points
of two elliptic normal curves $D_1,D_2 \in \mid 2H-R \mid$ on a cubic scroll 
$S \subset \Pn 4$, where $H,R \in PicS$ is the hyperplane class respectively
ruling on S. The resolution
$$0 \leftarrow {\cal O}_\Gamma \leftarrow {\cal O}_S 
\leftarrow 2{\cal O}_S(-2H+R) \leftarrow {\cal O}_S(-4H+2R) \leftarrow 0$$
of $\Gamma$ on $S$ induces a filtration

$$ \matrix{ 1 & - & - & - & - & \cr
            - & 3 & 2 & - & - & \cr
            - & - & - & - & - &  \cr
             } + 2 
\pmatrix{ - & - & - & - & - & \cr
            - & 2 & 3 & - & - & \cr
            - & - & - & 1 & - &  \cr 
             } +
\matrix{ - & - & - & - & - & \cr
            - & - & - & - & - & \cr
            - & - & 3 & 6 & 3 &  \cr
             }$$
for the syzygies of $\Gamma \subset \Pn 4$
$$ \matrix{ 1 & - & - & - & - & \cr
             - & 7 & 8 & - & - & \cr
            - & - & 3 & 8 & 3 & . \cr
             }$$
Hence the piece
$$ \matrix{ - & - & - & - & - & \cr
            - & - & - & - & - & \cr
            - & - & 3 & 8 & 3 & . \cr
             }$$
of the syzygies of $\Gamma$ recovers S as the support of the cokernel
of the map  $ 3{\cal O}(-4) \leftarrow 8{\cal O}(-5)$. The pencil can be 
recovered
from $2{\cal O}(-5)$ in the kernel of this map.\par
In the rest of this section we show that any general set of 8 points in $\Pn 4$ has syzygies as above. 

\vskip10pt 
\proclaim 5.2 Proposition.  There is a pencil of elliptic quintic curves through 8 points in
general position in $\Pn 4$.  This pencil sweep out a rational cubic scroll.\par
$Proof.$ Recall the
classical result (Castelnuovo) that there is a unique rational normal quartic curve through 7
general points in $\Pn 4$.  So given 8 points consider the rational normal curve passing through
all but one of them, say the point $p$.  The ideal of the rational quartic curve is generated by
the $2\times 2$ minors of a $2\times 4$ matrix with linear entries. This matrix has rank 2 at
$p$, but clearly one can find
 a $2\times 3$
submatrix which have rank 1 at $p$.
This submatrix has rank 1 on a rational cubic scroll. 
Thus we find rational cubic scrolls through the 8 points.\par

Given a rational cubic scroll through the 8 points there is either one or a pencil of
anticanonical curves on the scroll passing through the 8 points.
On the other hand, an anticanonical curve on the scroll is an elliptic quintic normal curve.  Any complete
linear series of degree 2 on this curve define a scroll, the union of the secant lines to the curve meeting
the curve in pairs of points which belong to the linear series. These scrolls are rational cubic scrolls, and
their union is precisely the secant variety of the elliptic curve.  Since the complete linear series of degree
2 on the elliptic curve $E$ is parameterized by $Jac_2 E$ which is again an elliptic curve, we have found for
each elliptic curve through the 8 points an elliptic pencil of rational
 cubic scrolls through the 8 points.\par

Given 8 points on an elliptic quintic curve, the quadrics through the 8 points will
define a linear series of degree 2 on the curve.  Consider the rational cubic scroll of the curve determined
by this linear series.  On it the elliptic curve will move in a pencil through the 8 points. Thus there is a
linear pencil of elliptic curves through the 8 points and they sweep out a rational cubic scroll.
 To show that there is only one such pencil, consider two elliptic quintic curves through the 8 points.  The 5
quadrics through one of the curves define a linear series of degree 2 residual to the 8 points on the other, so
there are 3 quadrics containing both curves.  Since their union has degree 10, the three quadrics cannot define
a complete intersection, so in fact they define a cubic scroll.
 Clearly it is the one described
above, so it is unique. $\slutt$\medskip
With (5.2) it is natural to see $VSP(F,8)$ for a cubic threefold $F$ as a variety of lines in the 
variety of apolar elliptic quintic curves to $F$.  We come back to this in section 8 after
introducing the spinor variety and the apolar Artinian Gorenstein ring 
$A^F$. 
\vskip20pt

{\bf 6. Equations and geometry of the spinor varieties ${\Sp_{ev}}$ and ${\Sp_{odd}}$}

{\vskip10pt \bf 6.0 } Let $W$ be a 10-dimensional vector space and let $Q^8$ be a rank 10 
quadric in
$\Pn 9={\bf P_{\ast}}(W)$, i. e. defined by a nondegenerate quadratic form $q:W\to k$. 
The quadric
$Q$ contains two families of isotropic $\Pn 4$'s,  parameterized by the (isomorphic) 
spinor varieties
${\Sp_{ev}}$ and ${\Sp_{odd}}$. We recall the description of these varieties cf. [Room 1952], [Chevalley
1954], [Lazarsfeld, Van de Ven 1984],[Zak 1984], [Ein 1986], [Mukai 1995].
 
{\vskip10pt \bf 6.1 }
Let $Q$ be defined by the equation
$$q=\sum_{i=1}^5 y_iy_{-i} = 0$$
Then $U_0 = \{y_1=\ldots y_5=0\}$ and $U_\infty = \{y_{-1}= \ldots y_{-5} = 0 \} \subset W$ represent a pair
of disjoint isotropic $\Pn 4$'s. Let $e_{\pm 1},\ldots,e_{\pm 5} \in W$ be the dual basis of the 
$y_{\pm i}$'s. Then $$U_0 = <e_{-1},\ldots,e_{-5}> \hbox{and } U_\infty = <e_1,\ldots,e_5>.$$
These space are dual to each other via q. Consider the Clifford-operators
$$\partial_+=\sum_{i=1}^5y_{i}e_i \hbox{ and } \partial_-=\sum_{i=1}^5y_{-i}e_{-i}.$$
 The operation $\partial_\pm \colon \wedge^{\bullet} U_\infty \to W^* \otimes 
\wedge^{\bullet \pm 1}U_\infty$
as wedge-product respectively contraction induce a map of bundles 
$\Lambda_{\Pn 9} \to \Lambda_{\Pn 9}(1)$, where $\Lambda = \wedge^{\bullet} U_\infty$.
For $\partial = \partial_+ + \partial_-$  one has
$\partial \circ \partial = q \cdot id_\Lambda$. Hence decomposing in even and odd parts
 one obtains
two rank 8 vector bundles on $Q$:
$$ E^{ev} = coker(\partial^{odd}\colon \Lambda^{odd}_Q(-1) \to \Lambda^{ev}_Q) \hbox{ and }
 E^{odd} = coker(\partial^{ev}\colon \Lambda^{ev}_Q(-1) \to \Lambda^{odd}_Q).$$
The bundles are related to each other via the short exact sequences
$$ 0 \to E^{odd}(-1) \to\Lambda^{ev}_Q \to E^{ev} \to 0 \hbox{ and } 
0 \to E^{ev} \to\Lambda^{odd}_Q(1) \to E^{odd}(1) \to 0. $$
$E^{ev/odd}$ are homogeneous bundles under the action of SO(q). They are the unique
indecomposable vector bundles $E$ on 
$Q$ with $H^i(Q,E(n))=0$  for all n and $0 < i < \hbox{dim }  Q$, cf. [Kn\"orrer 1987].
The induced action  on $\Lambda^{ev/odd} = H^0(Q,E^{ev/odd})$ are the half
spin representations of SO(q). 

{\vskip10pt \bf 6.2 } Every isotropic $\Pn 4$ arises as zero loci of a section in $E^{ev}$ 
respectively $E^{odd}$. We treat $E^{ev}$:

Let $1,e_1\wedge e_2,\ldots,e_4\wedge e_5, e_2 \wedge e_3 \wedge e_4 \wedge e_5 ,\ldots,
e_1 \wedge e_2 \wedge e_3 \wedge e_4$ be the basis of monomials in 
$e_1,\ldots e_5$ of $\Lambda^{ev}$, and let $x_0,x_{12},\ldots,x_{45},x_{2345},\ldots,x_{1234}$
be corresponding dual coordinate functions on $\Lambda^{ev}$. In these coordinates
the section $(1,0,\ldots,0)$ vanishes along $\P_*(U_0)$:
$$ \partial^{ev} \colon (1,0,\ldots,0) \mapsto y_1e_1+\ldots + y_5e_5.$$
Let $(id_5,A)=(\delta_{ij},a_{ij})$ be point in the affine neighborhood of
$0 \in Hom(U_0,U_\infty) \subset \Gr(5,W)$. $(id,A) \in Hilb_{\Pn 4} (Q)$ iff
$$\pmatrix {id, A\cr} \pmatrix{ 0 & id \cr id & 0 \cr} \pmatrix{id \cr A^t \cr} = A+A^t=0$$
i. e. if $A=(a_{ij})$ is skew symmetric. More intrinsically, an affine peace of
the spinor variety $\Sp_{ev}$ is given by $\Lambda^2U_\infty \subset
Hom(U_\infty^*,U_\infty) = Hom(U_0,U_\infty).$
Concerning the spinor embedding we consider the exponential map
$$exp: (a_{ij}) \mapsto (1,a_{ij},{\rm Pfaff}_{ijkl}(A)) \in \Lambda^{ev}.$$
Here ${\rm Pfaff}_{ijkl}(A)=a_{ij}a_{kl}-a_{ik}a_{jl}+a_{il}a_{jk}$. In more intrinsic terms
the exponential map is given by
$$exp \colon \Lambda^2U_\infty\to \Lambda^{ev}, A \mapsto 1+A+{1 \over 2}A \wedge A$$
A computation shows, that the corresponding section has zero loci defined by
$$y_i+\sum_j a_{ij}y_{-j} = 0 \hbox{ for } i=1,\ldots,5 $$
as desired.

{\vskip10pt \bf 6.3 } The image of the exponential map satisfies the homogeneous equations
$$\eqalign 
{
q^+_1&=x_0x_{2345}+x_{23}x_{45}-x_{24}x_{35}+x_{34}x_{25}\cr
q^+_2&=x_0x_{1345}-x_{13}x_{45}+x_{14}x_{35}-x_{34}x_{15}\cr
q^+_3&=x_0x_{1245}+x_{12}x_{45}-x_{14}x_{25}+x_{24}x_{15}\cr
q^+_4&=x_0x_{1235}-x_{12}x_{35}+x_{13}x_{25}-x_{23}x_{15}\cr
q^+_5&=x_0x_{1234}+x_{12}x_{34}-x_{13}x_{24}+x_{23}x_{14}\cr
q^-_1&=x_{12}x_{1345}+x_{13}x_{1245}+x_{14}x_{1235}+x_{15}x_{1234}\cr
q^-_2&=-x_{12}x_{2345}+x_{23}x_{1245}+x_{24}x_{1235}+x_{1234}x_{25}\cr
q^-_3&=-x_{13}x_{2345}-x_{23}x_{1345}+x_{34}x_{1235}+x_{1234}x_{35}\cr
q^-_4&=-x_{14}x_{2345}-x_{24}x_{1345}-x_{34}x_{1245}+x_{1234}x_{45}\cr
q^-_5&=-x_{15}x_{2345}-x_{25}x_{1345}-x_{35}x_{1245}-x_{1235}x_{45}\cr
}
$$
Indeed in the open set $x_0=1$ the first five equation express $x_{ijkl}$ as the Pfaffians of
$X=(x_{ij})$, and the last 5 equations become the well-known syzygies among the Pfaffians.

\proclaim 6.4 Proposition.\hbox{\rm (e.g. [Mukai 1995])}. $\Sp_{ev}\subset\Pn {15}$ is defined by the
equation above.  It has degree 12, its homogeneous coordinate ring is Gorenstein, with syzygies
$$ \matrix{ 1 & - & - & - & - & - \cr
         - & 10 & 16 & - & - & - \cr
            - & - & - & 16 & 10 & - \cr
              - & - & - & - & - & 1 \cr }$$\par
$Proof.$ We already saw that the spinor variety $\Sp_{ev} \subset \Pn {15}$ 
is defined by these equations in the affine part $x_0=1$. Moreover the equation have syzygies
as claimed. This follows from e.g. a MACAULAY computation. By the length of the resolution we conclude
that this ideal is projectively
Cohen-Macaulay. Furthermore by a similar computation $x_0$ is a nonzero divisor. 
Since projective Cohen-Macaulay ideals
are unmixed, this is the homogeneous ideal of $\Sp_{ev} \subset \Pn {15}$.$\slutt$ \par 

{\vskip10pt \bf 6.5 } Consider the point $o=(1:0:\ldots:0) \in \Sp_{ev}$ corresponding to $U_0$.
The tangent space at $o$ is $T_{\Sp,o} = \{x_{1234}=x_{1235}=x_{1245}=x_{1345}=x_{2345}=0 \}$. The equations

$$
\pmatrix {q^-_1\cr q^-_2\cr q^-_3\cr q^-_4\cr q^-_5\cr}= 
\pmatrix {0&x_{12}&x_{13}&x_{14}&x_{15}\cr
-x_{12}&0&x_{23}&x_{24}&x_{25}\cr
-x_{13}&-x_{23}&0&x_{34}&x_{35}\cr
-x_{14}&-x_{24}&-x_{34}&0&x_{45}\cr
-x_{15}&-x_{25}&-x_{35}&-x_{45}&0\cr} 
\cdot \pmatrix {x_{2345}\cr
x_{1345}\cr x_{1245}\cr x_{1235}\cr x_{1234}\cr }=0$$
define the image $\tilde \Sp \subset \Pn {14}$ of the birational projection of  $\Sp_{ev}$ from $o$. Notice that
$\tilde \Sp$ is generic syzygy variety, cf. [Eusen, Schreyer 1995]. $\tilde \Sp$ has syzygies
$$ \matrix{ 1 & - & - & - & - \cr
         - & 5 & 1 & - & - \cr
            - & - & 11 & 10 & 1  \cr
              - & - & - & - & 1  \cr }$$\par
The equation of $\Sp_{ev} \subset \Pn {15}$ can be recovered from $\tilde \Sp$ as
$$\varphi_4\cdot\pmatrix{x_0 \cr 1 \cr} = 0$$
where $\varphi_4$ is the last syzygy matrix of $\tilde \Sp$. 

The exceptional set
of the projection is inside $T_{\Sp,o}$ and there given by
$$\pmatrix {
x_{23}x_{45}-x_{24}x_{35}+x_{34}x_{25}\cr
-x_{13}x_{45}+x_{14}x_{35}-x_{34}x_{15}\cr
x_{12}x_{45}-x_{14}x_{25}+x_{24}x_{15}\cr
-x_{12}x_{35}+x_{13}x_{25}-x_{23}x_{15}\cr
x_{12}x_{34}-x_{13}x_{24}+x_{23}x_{14}\cr
}=0$$
These are just Pl\"ucker quadrics.
So the exceptional set is the seven dimensional cone
$$\Gr_o^7=\{\P_*(U)|{\rm dim} U\cap U_0\geq 3\}$$ 
consisting of even  isotropic $\Pn 4$'s whose intersection with 
$\P_*(U_0)$ is at
least 2-dimensional. Thus $\Gr_o^7$ has vertex $o$ and is isomorphic to the cone 
over $\Gr (3,U_0)$.

{\vskip10pt \bf 6.6 } The rational map 
$$ v^+:\P_*(\Lambda^{ev}) - - > \Pn 9=\P_*(W), a \mapsto [q_1^-(a):\ldots:q_5^+(a)]$$
defined by the 10 quadrics maps a section of $a \in \P_*(\Lambda^{ev}) =\P_*(H^0(Q, E^{ev}))$ to its zero loci,
which is a point on $Q$ for general sections since $c_8(E^{ev})=1$. Indeed
$$q_1^-\cdot q_1^++\ldots+q_5^-\cdot q_5^+=0,$$
so the image is $Q \subset \Pn 9$. To see the statement about the zero loci
we notice that from a different point of view the map
$$\partial^{ev} \colon \Lambda^{ev} \to W^* \otimes \Lambda^{odd}$$
coincides with the first syzygy matrix of $\Sp_{ev}$
$$ \varphi_2 \colon \Lambda^{ev} \to W^* \otimes (\Lambda^{ev})^*$$
if we identify  $\Lambda^{odd} \cong (\Lambda^{ev})^*$ via the wedge-product pairing into
$\Lambda^5U_\infty$. For a point $[a] \in \Pn {15} \backslash \Sp_{ev}$ the matrix $\varphi_2(a)$
has rank 9 precisely and $[q_1^+(a):\ldots:q_5^+(a)]$ represents the cokernel. This means that
the section $[a] \in \Pn {15}$ vanishes at $[q_1^+(a):\ldots:q_5^-(a)] \in Q$.

{\vskip10pt \bf 6.7 } Consider the inverse image by the map $v^+$ of an isotropic $\Pn 4$. In both
families this can be computed from the equations.
For the isotropic subspace $P_\infty=\P_*(U_\infty)\subset Q$ the inverse image by $v^+$  is defined by 
$$q^+_1=q^+_2=q^+_3=q^+_4=q^+_5=0$$
Let $H$ be the hyperplane defined by $x_0=0$.  
The quadrics, restricted to $H$, are the same Pl\"ucker quadrics as above.  
This time they define an
11-dimensional cone $G_{P_\infty}^{11}$ isomorphic to a cone over $\Gr (2,U_\infty)$.
Note that the hyperplane $H$ belongs to the dual variety.  It is tangent to ${\Sp_{ev}}$ along
 the contact
locus $C_H$, the $\Pn 4$ which is the vertex of the 11-dimensional cone $G_{P_\infty}^{11}$, it is 
defined by $x_{ij}=0, {\rm all}\quad i<j$ and coincides naturally with the dual space
to $\P_*(U_\infty)$. On the complement of $H$, where $x_0=1$, we get as noted in $(6.3)$
precisely the image of the exponential map. From the equations we get that the quadrics
$q^+_i$ define the graph in ${\rm span} <T_{\Sp ,o} ,C_H>=\Pn {15}$ of the rational map of 
$T_{\Sp ,o}$
into $C_H$ defined by the Pl\"ucker quadrics. This graph is 10 dimensional and its closure in
$\Pn {15}$ is of course again just ${\Sp_{ev}}$ (cf. also [Ein 1986]). Thus the inverse image of
$P_\infty$ is the 11-dimensional cone $G_{P_\infty}^{11}$ which intersect ${\Sp_{ev}}$ along a tangent
hyperplane section $H$ defined by $x_0=0$.
Therefore the isotropic space ${P_\infty}$ canonically correspond to the hyperplane
which contains its inverse image by $v^+$. The point $H\in {\Sp_{odd}}$ is the point
 associated to $P_\infty$.

{\vskip10pt \bf 6.8 } For the isotropic space $P_0={\bf P_{\ast}}(U_0)$, the inverse image
 by $v^+$ is given by
$$q^-_1=q^-_2=q^-_3=q^-_4=q^-_5=0. $$
As noted in (6.5) this is the cone with vertex $o=(1:0:\dots :0)\in {\Sp_{ev}}$ over a generic syzygy
variety.  
Two isotropic subspaces of $Q$ which intersect in a $\Pn 3$ form the intersection of $Q$ with a $\Pn 5$.
Therefore the union of the preimages of the two subspaces by $v^+$ is the complete intersection of 4
quadrics. Since the preimage of $P_\infty$ has degree 5 and dimension 11 the preimage of $P_0$ must
therefore have degree 11 and dimension 11. 
  
Note that any hyperplane containing $T_{\Sp ,o}$ will intersect the preimage of $P_0$ in $T_{\Sp ,o}$ and a residual variety which
intersect $T_{\Sp ,o}$ along a quadric.
\vskip20pt
{\bf 7. The apolar Artinian Gorenstein ring of a general cubic threefold}
\bigskip

 {\it Proof of }(0.10):
Recall from $(6.4)$ that the homogeneous coordinate ring $R_S$ of the spinor variety 
$\Sp\subset\Pn {15}$ is Gorenstein with syzygies
$$ \matrix{ 1 & - & - & - & - & -& \cr
           - & 10 & 16 & - & - & -& \cr
            - & - & - & 16 & 10 & -& \cr
              - & - & - & - & - & 1&. \cr }$$
The spinor variety has dimension $10$, so eleven general linear forms $h_0,\dots,h_{10}$ define a
$\Pn 4$ which does not intersect $\Sp$.  Therefore the quotient $A=R_S/(h_0,\dots,h_{10})$ is an
Artinian Gorenstein ring.  Its Hilbert function is $(1,5,5,1)$ and it has socle degree 3, so $A$
is the apolar Artinian Gorenstein ring $A^F$ for some cubic hypersurface $F\subset \Pnd 4$.
Now, let $X$ be the open set in the Grassmannian $\Gr (5,\Lambda^{ev})$ of subspaces which
does not intersect the spinor variety $\Sp+$.  The action of $SO(W,q)$ on $\Lambda^{ev}$ induces an
action on $X$.  Clearly, any isomorphism between two quotients $A$ and $A^{\prime}$ is induced by an
automorphism of $\Pnd 4$. Therefore, the above construction defines a map  $m:X/SO(W,q)\to H_3$, to
the space of cubics in $\Pn 4$ modulo projective equivalence. That this map is injective follows
from the
\bigskip \proclaim 7.1 Lemma.  For any cubic $F$ in the image of $m$ there is precisely one
quadratic relation between the quadrics in $F^{\bot}_2$\par
The proof of this lemma is postponed until (8.7).\medskip
   The dimension of $X/SO(W,q)$ is
dim$\Gr(5,16)-{\rm dim}SO(W,q)=55-45=10$ while dim$H_3=\coh 0{\Pn 4}3-{\rm dim}GL(\C^5)=35-25=10$ so
the map $m$ is dominating onto its image components.  But both source and target are irreducible, so
the image of $m$ must be a dense subvariety, and $(0.10)$ follows.$\slutt$
\bigskip
We conjecture that the generality condition on $A$ can be made
precise via syzygies:
\bigskip \proclaim 7.2 Conjecture. An Artinian Gorenstein ring A with Hilbert function
$(1,5,5,1)$ is isomorphic to a quotient $R_S/(h_0,\ldots,h_{10})$ for some
 linear forms
$h_0,\ldots,h_{10} \in R_\Sp$ iff $A$ has syzygies
$$ \matrix{ 1 & - & - & - & - & -& \cr
           - & 10 & 16 & - & - & -& \cr
            - & - & - & 16 & 10 & -& \cr
              - & - & - & - & - & 1&. \cr }$$
as $\C[x_0,\ldots,x_4]$-module. \par

This conjecture implies  Mukai's Theorem (Theorem 0.9) by applying the
concept of `very rigidity' cf. [Buchweitz 1981].
For a different proof of Mukai's Theorem see also [Eusen,Schreyer 1995].

\vskip20pt {\bf 8.  Proof of the main results }
\bigskip

{\bf 8.0 }  Let $A^F$ be the
homogeneous
coordinate ring of the {\sl empty} intersection of the spinor 
variety $\Sp$ with a general $P=\Pn 4$, i. e. the apolar Artinian Gorenstein
 ring of a general cubic $F\in \Pn 4$.

To describe $VSP(F,8)$ we follow the procedure of section 3 and consider duality and the 
projection from
partials $\pi_2^F$.   The projection from partials $\pi_2^F:P\to\Pn 9$ 
is nothing but
the restriction of the map $v^+$ to $P$, and the points in  $VSP(F,8)$
representing smooth apolar subschemes to $F$ of degree 8 correspond to proper 8-secant planes to
the image $V_P=\pi^F_2(P)\subset \Pn 9$. We show in lemma 8.4 below that any such 8-secant plane is 
contained in the quadric $Q\subset \Pn 9$ which is the image of $v^+$.  Since any plane in $Q$ is the
intersection of  a pencil of isotropic
 $\Pn 4$'s, the preimage of a plane is contained in a pencil of tangent hyperplanes to $\Sp$, cf.
6.7.  Thus
any plane in $Q$ correspond to a line in the dual variety ${\Sp_{odd}}$ to $\Sp={\Sp_{ev}}$.\par
Recall that the preimage of an isotropic $\Pn 4$ corresponding to a point on ${\Sp_{odd}}$ is a
11-dimensional cone over the Grassmannian variety $\Gr (2,5)$ contained in a tangent hyperplane. 
Clearly $P$ intersects this inverse image properly in a curve precisely when $P$ is contained in
the  tangent hyperplane.   In this case the curve of intersection is an elliptic quintic
curve.  Therefore   $Y=\Pn
{10}\cap{\Sp_{odd}}\subset{\bf P_{\ast}}(\Lambda^{odd})$  where $\Pn {10}$ is dual to $P\subset {\bf
P_{\ast}}(\Lambda^{ev})$, parameterizes apolar elliptic quintic curves to $F$, i. e. defines a
subvariety of $VPS(F,5t)$ in the notation of 1.3. \par On the other hand, it follows from 
(5.2) that apolar subschemes of degree 8 are contained in a pencil of apolar elliptic quintic curves. 
Let $M_Y$ denote the Fano variety of lines in $Y$.  We show that for general cubic $3$-fold $F$, the
variety $M_Y$ is smooth of  dimension $5$ and that every point on $M_Y$ correspond to a finite apolar
subscheme to $F$ of length 8, i. e. that $F^5$ coincide with $VSP(F,8)$. \par

\proclaim 8.1 Lemma. Let $Y=\Pn {10}\cap {\Sp_{odd}}\subset{\bf P_{\ast}}(\Lambda^{odd})$ 
for a general $\Pn {10}$, and
let $M_Y$ be the Fano variety of lines in $Y$.  Then $M_Y$ is smooth of dimension $5$.\par
$Proof.$  Recall for each point $H\in {\Sp_{odd}}$ the tangent cone $G_H^7= {\Sp_{odd}}\cap T_H$  
with vertex $H$ inside the tangent space
$T_H$ is isomorphic to a 7-dimensional cone over the Grassmannian of planes in the isotropic
$\Pn 4$ corresponding to $H$. For each point $H$ in $Y$ let $g_H^2=G_H^7\cap \Pn {10}\subset Y$.
Then  $g_H^2$ is the cone over an elliptic normal quintic curve.  In particular there 
is an elliptic
curve of lines through a general point $H$ in $Y$. Since a line is a pencil of points, 
$Y$ contains
a  5-dimensional family of lines, i. e. $M_Y$ has dimension $5$.\par
For smoothness we consider the normal bundle $N_L$ of a line $L$ in $Y$.  $M_Y$ is 
smooth as soon as $h^1(N_L)=0$
for every line $L$ in $Y$.  
This is verified by a dimension count. For each line $L$ in $\Sp$ the normal bundle is $3{\cal O}_L
+ 6\O_L(1)$ (cf. [Ein 1986]).  Since $Y$ is a section of ${\Sp_{odd}}$ with a $\Pn {10}$, the normal bundle
sequence becomes
$$ 0\longrightarrow \; N_L\; \longrightarrow \; 3{\cal O}_L\oplus 6{\cal O}_L(1)\;\longrightarrow \;
5{\cal O}_L(1) \; \longrightarrow 0 $$
This is exact so the rank of the map $\alpha_L :6{\cal O}_L(1)\;\longrightarrow \; 5{\cal O}_L(1)$
is at least 3;  otherwise the map would have a summand $3{\cal O}_L \;\longrightarrow \; 3{\cal
O}_L(1)$ which cannot be surjective at every point of $L$.  On the other hand $h^1(N_L)\not= 0$ 
only if $N_L$ has
a summand ${{\cal O}_{L}(-d)}$ with $d\geq 2$.  Therefore  $L$ is a singular point 
on $M_Y$ only if the map
$\alpha_L$ has rank 3. But the set of flags $(L,\Pn {10})$ of lines $L$ in 
${\Sp_{odd}}\cap\Pn {10}$ where $\alpha_L$
has rank 3 is a determinantal variety of codimension 6 in the variety of all
 flags so for general
$Y$, $M_Y$ is smooth.$\slutt$\par
\medskip

For the correspondence between $8$-secant planes to $V_P$ and lines in $Y$, we first study 
the fibers of the
map $v^+$ more carefully.\par
 Let $Bl_{\Sp}:\tilde {\bf P }\to {\bf P_{\ast}}(\Lambda^{ev})$ be the the blowup map
along $\Sp$, then $v^+$ lifts to a morphism $\tilde {v^+}:\tilde {\bf P }\to Q\subset\Pn 9$.
  For a subscheme
$Z\subset Q$ let $F_Z=Bl_{\Sp}((\tilde {v^+})^{-1}(Z))$.  By
abuse of notation we often call $F_Z$ the preimage of $Z$ by $v^+$.  
We consider first the preimage of points,
lines and planes.

\proclaim 8.2 Lemma. The fiber of $v^+$, i. e. $F_p$ for a point $p\in Q$, is a linear
$\Pn 7$ which
intersect $\Sp$ in a smooth quadric hypersurface.  The preimage of a line, i. e. $F_L$
for a line
$L\subset Q$ is a  rational normal quartic scroll in a $\Pn {11}$ with vertex a $\Pn 3$. The
preimage of a conic section, which is not contained in $Q$, is a rational normal scroll of
dimension 8 and degree 8.\par $Proof$. By homogeneity we may choose any point and any line and compute
the preimage from the equations.\par Consider the point on $Q$ with $y_1=q^+_1$ as the only nonzero
coordinate. The other quadrics all vanish on the $\Pn
7=\{x_{12}=x_{1345}=x_{13}=x_{1245}=x_{14}=x_{1235}=x_{15}=x_{1234}=0\}$ which is
the singular locus of $q^-_1$ , and $q^+_1$ defines in it a smooth quadric. For any point outside
the $\Pn 7$ one of the variables defining it is nonzero.  But on the complement of any coordinate
hyperplane $\Sp$ is defined by 5 quadrics as explained in (6.3), and $q^+_1$ is not among
these 5 quadrics for any of the coordinates defining the $\Pn 7$.  Therefore any point outside this
$\Pn 7$ on all but the first quadric also lie on the first quadric. Thus we have recovered $F_p$ as
a $\Pn 7$ which intersect ${\Sp_{ev}}$ in a quadric.  In
 particular the fiber over a point intersect ${\Sp_{ev}}$ in codimension one.\par
Next look at the line in $Q$ with all coordinates but $y_1=q^-_1$ and $y_2=q^-_2$ equal zero. 
In this case the remaining quadrics  
$$\eqalign 
{
q^+_1&=x_0x_{2345}+x_{23}x_{45}-x_{24}x_{35}+x_{34}x_{25}\cr
q^+_2&=x_0x_{1345}-x_{13}x_{45}+x_{14}x_{35}-x_{34}x_{15}\cr
q^+_3&=x_0x_{1245}+x_{12}x_{45}-x_{14}x_{25}+x_{24}x_{15}\cr
q^-_3&=-x_{13}x_{2345}-x_{23}x_{1345}+x_{34}x_{1235}+x_{1234}x_{35}\cr
q^+_4&=x_0x_{1235}-x_{12}x_{35}+x_{13}x_{25}-x_{23}x_{15}\cr
q^-_4&=-x_{14}x_{2345}-x_{24}x_{1345}-x_{34}x_{1245}+x_{1234}x_{45}\cr
q^+_5&=x_0x_{1234}+x_{12}x_{34}-x_{13}x_{24}+x_{23}x_{14}\cr
q^-_5&=-x_{15}x_{2345}-x_{25}x_{1345}-x_{35}x_{1245}-x_{1235}x_{45}\cr
}
$$
cut the union of $\Sp$ and the rational quartic scroll in $\Pn {11}=\{x_0=x_{34}=x_{45}=x_{35}=0\}$ 
defined by the $2\times
2$ minors of: $$\pmatrix 
{
x_{13}&x_{14}&x_{15}&x_{1345}\cr
x_{23}&x_{24}&x_{25}&-x_{2345}.\cr
}
$$
In fact the 5 quadrics defining $\Sp$ on the complement of any of the coordinate hyperplanes which contain $\Pn
{11}$ does not involve $q^-_1$ or $q^-_2$, so any point outside the $\Pn {11}$ on the 8 quadrics lie on $\Sp$. 
  
Thus the preimage $F_L$ of the line is a rational normal quartic scroll in $\Pn {11}$ with vertex a $\Pn 3$.
 One can check that the
intersection of $\Pn {11}$ with $\Sp$ is not proper,
in fact this intersection has codimension one in the fiber
over the line.\par
For the preimage of a conic section, the computation is similar.
$\slutt$\medskip

A plane $\pi\subset Q$ is contained in a pencil of isotropic $\Pn 4$'s of each family on $Q$.  These
pencils of isotropic $\Pn 4$'s correspond to a line $L_{\pi}\subset \Sp$ and by duality a pencil of
tangent hyperplanes defining a $\Pn {13}_{\pi}\subset {\bf P_{\ast}}(\Lambda^{ev})$.  The
following lemma should be compared with (5.2).  

\proclaim 8.3 Lemma.   The preimage by $v^+$ of a plane, i. e. $F_{\pi}$ for a plane $\pi\subset Q$
is a cone $C_\pi$ over a codimension 4 subvariety of degree 8 and sectional genus 3 with vertex
line $L_\pi\subset \Sp$
 inside ${\bf P}_{\pi}^{13}$. More precisely, the tangent spaces to $\Sp $ along $L_{\pi}$ form the
rulings of a rational cubic croll and $F_{\pi}$ is the baselocus of a pencil of divisors
 on this scroll determined by the tangent hyperplanes $H\supset\Pn {13}_\pi$.  These divisors are
all isomorphic to the $10$-dimensional cone over a tangent hyperplane section of $\Gr (2,5)$. \par

$Proof.$ We consider the plane $\pi$ in $Q$ with all coordinates but
$y_1=q^-_1$, $y_2=q^-_2$ and $y_3=q^-_3$ equal zero.
  The remaining quadrics   $$\eqalign  {
q^+_1&=x_0x_{2345}+x_{23}x_{45}-x_{24}x_{35}+x_{34}x_{25}\cr
q^+_2&=x_0x_{1345}-x_{13}x_{45}+x_{14}x_{35}-x_{34}x_{15}\cr
q^+_3&=x_0x_{1245}+x_{12}x_{45}-x_{14}x_{25}+x_{24}x_{15}\cr
q^+_4&=x_0x_{1235}-x_{12}x_{35}+x_{13}x_{25}-x_{23}x_{15}\cr
q^-_4&=-x_{14}x_{2345}-x_{24}x_{1345}-x_{34}x_{1245}+x_{1234}x_{45}\cr
q^+_5&=x_0x_{1234}+x_{12}x_{34}-x_{13}x_{24}+x_{23}x_{14}\cr
q^-_5&=-x_{15}x_{2345}-x_{25}x_{1345}-x_{35}x_{1245}-x_{1235}x_{45}\cr
}
$$
cut the union of $\Sp$ and a cone $C_\pi$ over a codimension 4 subvariety of degree 8 and sectional
genus 3 with a vertex line $$L_\pi :\quad x_{ij}=0\quad {ij}\not= {45}\quad x_{2345}=x_{1345}=x_{1245}=0$$
in ${\bf P}_{\pi}^{13}=\{x_0=x_{45}=0\}$.  In fact any point outside ${\bf P}_{\pi}^{13}$ on the 7 quadrics lie on $\Sp$
like above.\par
In each tangent hyperplane $H\supset{\bf P}_{\pi}^{13}$ the preimage $G_H^{11}$ of the
isotropic $\Pn 4$ corresponding to $H$ is defined by 5 quadrics.  For the pencil of tangent
hyperplanes there is a net of common quadrics, the $2\times 2$ minors of the matrix
$$\pmatrix 
{
x_{14}&x_{24}&x_{34}\cr
x_{15}&x_{25}&x_{35}\cr
}
$$ which define a rational cubic scroll of dimension 11 with vertex a $\Pn 7$ inside $\Pn
{13}_{\pi}$. For each $H$, the intersection  $G_H^{11}\cap{\bf P}_{\pi}^{13}$ has dimension 10 and is a
divisor on this cubic scroll. Now $G_H^{11}$ is naturally the cone over the Grassmannian $\Gr (2,5)$
 of lines in the isotropic $\Pn 4$ corresponding to $H$, while the pencil of isotropic $\Pn 4$'s
corresponding to the pencil of hyperplanes through ${\bf P}_{\pi}^{13}$ meet in a plane.  Therefore
$G_H^{11}\cap{\bf P}_{\pi}^{13}$ is isomorphic to the cone over a tangent
hyperplane section of $\Gr (2,5)$.
 The preimage $F_{\pi}$ 
  is the intersection on cubic scroll of this pencil of divisors.\par
Explicitly $F_\pi$ is defined by the above net of quadrics and the 4 quadric entries of the product
matrix : $$\pmatrix 
{
x_{14}&x_{24}&x_{34}\cr
x_{15}&x_{25}&x_{35}\cr
}\cdot \pmatrix 
{
x_{23}&x_{2345}\cr
-x_{13}&x_{1345}\cr
x_{12}&x_{1245}\cr
}.
$$
 Each of the divisors on the cubic scroll is defined by an element in the column space of the second factor
above.  The vertex of this divisor as a cone is the $\Pn 4$
defined by the corresponding column vector in the vertex $\Pn 7$ of the cubic scroll. \par 
 A point on the line $L_\pi$ correspond to a vector in the row space of the first matrix
factor above.  In fact the entries in the row vector vanish in ${\bf P}_{\pi}^{13}$ exactly along the
tangent $\Pn {10}$ to $\Sp$ at the corresponding point. The inverse image of the isotropic subspace
with  vertex at a point
on the line
$L_\pi$ restricts to the ${\bf P}_{\pi}^{13}$ as the union of this tangent
space and the cone $F_\pi$.  In fact the tangent spaces to $\Sp$ along
the vertex line form the ruling of the distinguished cubic scroll that contains $F_\pi$, 
and the $\Pn 7$ inside $F_\pi$
is precisely the intersection of these tangent spaces.$\slutt$\medskip
  Thus $v^+$ identifies the line in both ${\Sp_{ev}}$ and ${\Sp_{odd}}$ which correspond to a plane in
$Q$. One can check that the intersection of ${\Sp_{ev}}$ with the inverse image $F_{\pi}$ has codimension
2. \par 
Next we restrict $v^+$ to the $P=\Pn 4\subset {\bf P_{\ast}}(\Lambda^{ev})$ which
correspond to the apolar
Artinian Gorenstein ring $A^F$.  Note that $P$ does not intersect $\Sp$, and $v^+$ have only
linear fibers which
intersect $\Sp$ in codimension 1, so this restriction is an embedding.  Like above we set $V_P=v^+(P)\subset
Q\subset \Pn 9$.
\proclaim 8.4 Lemma.  Any proper 8-secant plane to $V_P$ is contained in the quadric $Q$.
\par
$Proof$. Since $V_P$ is contained in $Q$, the intersection of $V_P$ with any plane which is not
contained in $Q$ is contained in a conic section. But the preimage of a conic section is a rational
normal scroll of degree 8 and codimension 7 by (8.2). This scroll cannot meet any $\Pn 4$ in 8
points, so the lemma follows.$\slutt$

\proclaim 8.5 Lemma. Assume that $P$ does not intersect the spinor variety, and that the
intersection $Y$ of the
dual $\Pn {10}=P^{\bot}\subset {\bf P_{\ast}}(\Lambda^{odd})$ with the dual spinor variety is smooth
and has a
smooth Fano variety of lines $M_Y$.  Then every line in $Y$ correspond to a plane in $Q\subset\Pn 9$
which
intersect $V_P$ along a finite scheme of length 8.\par

$Proof.$ 
Any line $l_\pi\subset Y$ correspond to a dual $\Pn {13}\subset{\bf
P_{\ast}}(\Lambda^{ev})$ which contain $P$.  The line $l_\pi$ correspond to a plane 
$\pi\subset Q\subset \Pn 9$ whose
preimage $F_\pi$ by $v^+$ span this $\Pn {13}$.  The preimage of $\pi\cap V_P$
is therefore the
intersection $P\cap F_{\pi}$.  But $F_\pi$ has degree 8 and
codimension 4, by (8.3), so as soon as the intersection $P\cap F_{\pi}$ is proper the lemma
follows.\par

On the other hand since $v^+$ restricted to $P$ is an embedding, the intersection $\pi\cap V_P$
contains at most
a curve.  And, since the map is defined by quadrics, any plane curve on $V_P$ is a conic section, the
image of a line on $P$.  Therefore we need only to rule out the existence of a line $L$ in the intersection
 $P\cap F_{\pi}$.\par
Assume that $P\cap F_{\pi}$ contains a line $L$.
For a case by case argument we use the coordinates of the proof of (8.3).  Thus $\Pn {13}$ is defined
by  $\Pn {13}=\{x_0=x_{45}=0\}$ and $F_\pi$ is defined by the by the $2\times 2$ minors of the matrix
$${\cal M}=\pmatrix 
{
x_{14}&x_{24}&x_{34}\cr
x_{15}&x_{25}&x_{35}\cr
}$$ 

and the quadratic entries of the product matrix:
$${\cal M}\cdot {\cal N}=\pmatrix 
{
x_{14}&x_{24}&x_{34}\cr
x_{15}&x_{25}&x_{35}\cr
}\cdot \pmatrix 
{
x_{23}&x_{2345}\cr
-x_{13}&x_{1345}\cr
x_{12}&x_{1245}\cr
}.
$$
The intersection of $F_\pi$ with $\Sp$ is defined with the additional three quadrics
$$
\eqalign { 
q^-_1&=x_{12}x_{1345}+x_{13}x_{1245}+x_{14}x_{1235}+x_{15}x_{1234}\cr
q^-_2&=-x_{12}x_{2345}+x_{23}x_{1245}+x_{24}x_{1235}+x_{1234}x_{25}\cr
q^-_3&=-x_{13}x_{2345}-x_{23}x_{1345}+x_{34}x_{1235}+x_{1234}x_{35}\cr
}
$$
The $2\times 2$-minors of ${\cal M}$ define the distinguished cubic scroll $S_\pi$ which contain $F_\pi$.

The row vectors of ${\cal M}$ define the ruling of the cubic scroll and it is the
intersection with tangent spaces to ${\Sp_{ev}}$ with the $\Pn {13}$.  
Since any tangent space intersect $\Sp$ in codimension 3, the intersection of $P$ with a tangent
space is at most a plane, i. e. the rank of any row vector of ${\cal M}$ restricted to $P$ is at least
2. If every row vector has rank 2 restricted to $P$, then $S=S_\pi\cap P$ is an irreducible quadric
hypersurface in $P$, otherwise $S$ is a possibly reducible cubic surface scroll.
 \par
 Now, the line $L$ is clearly contained in $S$. Consider the rank of a general row vector of ${\cal M}$ restricted to
$L$.  This rank may be 0,1 or 2.  If the rank is 0 then $S$ is a cone with $L$ inside the vertex.  If the rank
is 1 then $L$ intersect the general ruling of $S$ in a point, and if the
rank is 2 then $L$ is contained in one ruling and do not intersect any other ruling of $S$. We go
case by case.

When ${\cal M}$ have rank 0 on $L$, then the three additional quadrics reduce to the $2\times 2$
minors of ${\cal N}$.  Therefore $L$ need not intersect the spinor variety.  But we will show
that the line $l_\pi\subset Y$ in this case correspond to a singular point on the Fano variety
$M_Y$ of lines on $Y$.
  Now, the rulings in $S$ all pass through $L$; if $S$ is a surface it is three plans through $L$, if it is a
quadric it is a cone with vertex $L$, in either case it contains three planes which together span $P$.  Each
plane is the intersection of $P$ with a tangent space, so dually each of these planes correspond to planes in
$Y$ through the line $l_\pi$. Moreover, these planes in $Y$ must span a $\Pn 4$; otherwise $Y$ contains a $\Pn
3$.  Any $\Pn 3\subset Y$ correspond by duality to a $\Pn {11}$ containing $P$ and the tangent space $T_p$ to
$\Sp$ at a point $p$.  Thus $P$ would intersect this tangent space in a $\Pn 3$.  But the tangent cone to
$\Sp$ has codimension 3 inside $T_p$, cf. (6.5), so $P$ would in fact intersect $\Sp$, which is
absurd.  Since the planes in $Y$ which contain $l_\pi$ span a $\Pn 4$, the normal bundle of $l_\pi$
in $Y$ has at least 6 sections, 2 from each of three independent planes.  But the Fano variety $M_Y$
has dimension 5 so the line $l_\pi$ must be a singular point on $M_Y$, contrary to our assumption.
\par 
 When the general row of ${\cal M}$ have rank 1
 on $L$, then we may assume that  $x_{14}$ and $x_{15}$ are the only nonzero terms of ${\cal M}$ on $L$.
The remaining quadrics defining $F_\pi$ restrict to $L$ as the quadratic forms  
$x_{14}x_{23}=x_{14}x_{2345}=x_{15}x_{23}=x_{15}x_{2345}=0$. The last three
quadrics defining the intersection with the spinor variety reduce to the quadric
$x_{12}x_{1345}+x_{13}x_{1245}+x_{14}x_{1235}+x_{15}x_{1234}$, so the line must intersect
the spinor variety, contrary to our assumption that $P$ does not.\par
  When the general row of ${\cal M}$ have rank 2
 on $L$, then we may assume that $x_{14}$ and $x_{24}$ are the only nonzero
terms on $L$.  The remaining quadrics defining $F_\pi$ restrict to $L$ as the quadric forms 
$x_{14}x_{23}-x_{24}x_{13}=x_{14}x_{2345}+x_{24}x_{1345}=0$. The last three quadrics defining
the intersection with the spinor variety reduce to the quadrics
$$x_{12}x_{1345}+x_{13}x_{1245}+x_{14}x_{1235}= -x_{12}x_{2345}+x_{23}x_{1245}+x_{24}x_{1235}=
 -x_{13}x_{2345}-x_{23}x_{1345}.$$ 
These  5 quadrics are the principal Pfaffians of the following
$5\times 5$-dimensional skew symmetric matrix: $$\pmatrix {0&0&-x_{24}&-x_{2345}&x_{23}\cr
0&0&x_{14}&-x_{1345}&-x_{13}\cr
x_{24}&-x_{14}&0&x_{1245}&-x_{12}\cr
x_{2345}&x_{1345}&-x_{1245}&0&-x_{1235}\cr
x_{23}&x_{13}&x_{12}&x_{1235}&0\cr}$$

The two quadrics defining $F_\pi$ are the last two principal 
Pfaffians of this matrix. When these vanish on the line $L$, then all the Pfaffians vanish 
on the point defined by  $x_{1235}$ on $L$, so again $L$ must intersect the spinor variety.\par
 We have shown that the intersection
$P\cap F_\pi\subset \Pn {13}$ is proper for any plane $\pi$ corresponding to a line in $Y$, so
 the lemma follows.$\slutt$\medskip

From (0.10), (8.2), (8.4) and (8.5) we get

\proclaim 8.6 Theorem. Let $F \subset \Pn 4$ be a general cubic. There exists a linear 
subspace $\Pn {10} \subset \Pn {15}$, which depends on $F$, such that $VSP(F,8)$ is isomorphic to 
the variety of lines  in the 5-fold $Y=Y(F) := \Pn {10} \cap \Sp^{10} \subset \Pn {15}$. 
Moreover $VSP(F,8)$ is smooth of dimension 5.\par\medskip
{\bf 8.7} {\it Proof of} (7.1). 
 For this proof we investigate the geometry of $M_Y$ a bit further.  Recall from the proof of (8.1)
 that every point $H$ in $Y$ is the vertex of a cone $g_H$ over an elliptic quintic curve.  But each
line in the cone $g_H$ correspond to a plane in $P_H$, the isotropic $\Pn 4$ corresponding to $H$. 
  So for each $H\in Y$ there is an elliptic quintic scroll $V_H$ of planes in $P_H$. 
This threefold scroll is singular along an elliptic quintic surface scroll $S_H$ of lines in $P_H$.  
Turning to the dual space, note that for a general $\Pn {10}$ the dual $P=\Pn 4$ do not intersect
${{\bf S}_{ev}}$ at all. The space $P$ meets each cone $G_H^{11}$ for a general $H\in {{\bf S}_{odd}}$
in 5 points, but for a point $H$ in $Y$, the span of $G_H^{11}$ contains $P$ and the intersection is
an elliptic curve $E_H$, the intersection of a Grassmannian with a general linear space. For two
points $H$ and $H^{\prime}$ on a line, $G_H^{11}$ and $G_{H^{\prime}}^{11}$ have a common line in the
vertex, it is the dual of the plane of intersection of $P_H$ and $P_{H^{\prime}}$.  Since $G_H^{11}$
is the cone over the Grassmannian of lines in $P_H$ each point on $E_H$ determines a line in
$P_H\subset Q $.  But $P$ lies in the span of both $G_H^{11}$ and $G_{H^{\prime}}^{11}$, so $P$
intersects $G_H^{11}$ only in fibers of lines in $P_H$ which meet $P_{H^{\prime}}$. The elliptic
curve $E_H$ altogether determines an elliptic scroll of lines. Since each line in the scroll meet the
plane $P_H\cap P_{H^{\prime}}$ for any $H^{\prime}\in g_H$, the scroll must coincide with the scroll
$S_H$ in $P_H$ described above.\par
Now, any line in $Y$ correspond to an $8$-secant plane to the image
$v^+(P)=V_P\subset \Pn 9$. On the scroll $S_H$ these planes are the planes of plane cubic curves. The
image of $E_H$ is an elliptic curve of degree 10 on the scroll $S_H$, which must therefore intersect
the plane cubic curves on $S_H$ in 8 points.  More precisely the image of $E_H$ is a curve of 
type $C_0+7R$ on $S_H$, where $C_0$ is a section with selfintersection 1 and $R$ is a member of the
ruling.  \par
 Since the scroll $S_H$ is not contained in any quadric, it follows that
the image of $E_H$ is not contained in any quadrics either.  Therefore the isotropic $\Pn 4$ of any
apolar elliptic curve $E_H$ is contained in any quadric which contains $V_P$.  But for a general
point on the quadric $Q$ there is a 6-dimensional quadric of isotropic $\Pn 4$'s through the point,
they form a 6-dimensional quadric in the spinor variety, and necessarily intersects the 5-dimensional
subvariety $Y$.  Therefore the isotropic subspaces of apolar elliptic curves fill all of $Q$, and
this is the unique quadric containing $V_P$.$\slutt$

\vskip20pt

{\bf 9. Invariants of $VSP(F,8)$}

{\vskip10pt \bf 9.1 } To find the invariants of $M_Y$ we describe the Fano variety $M$ of
lines in $\Sp (={\Sp_{ev}})$
via the correspondence with the Fano variety of planes in the quadric $Q$.
These varieties are isomorphic.  Each line in $\Sp$ correspond to a pencil of isotropic $\Pn 4$'s of
the same family through a plane in the quadric $Q$, and this correspondence is made precise by the
map $v^+$ (cf. 8.3).  In fact $v^+$ induces a map of $\Gr (3,{\cal N}(-2h)=\Gr (3,B)$ into $\Gr
(3,W)$, so that the universal subbundle on $\Gr (3,B)$ is induced from $M\subset \Gr (3,10)$. 
The Fano variety $F_q$ of planes in $Q$ is easily described in the Grassmannian $G=\Gr (3,10)$
of planes in $\Pn 9$.  But we need to identify the subvariety $M_Y$ so we compute $M\subset\Gr
(2,16)$.\bigskip

\proclaim 9.2 Proposition.  $M_Y$ is a Fano 5-fold of index 1 and degree 660.\par
$Proof.$ As above, let $M$ denote the Fano variety of lines in the spinor variety $\Sp$ and
let $G \subset \Sp \times M$ denote the universal family.
 $pr_1 \colon G \to \Sp$ is a $\Gr  (3,5)$-bundle. We compute the cohomology rings of $\Sp$, $G$ and
$M$,  starting with $\Sp$.

Let $B$ denote the tautological rank 5 subbundle on the spinor variety $\Sp$ with tautological
sequence 
$$ 0 \to B \to 10\O_\Sp \to B^* \to 0.$$

Then $\Omega^1_{\Sp} \cong \Lambda^2 B \subset {\cal H}om(B,B^*)^*$.
The cohomology ring $H^*(\Sp,\Q)$ is generated by the Chern classes of $B$ and the hyperplane
class $h$ of $\Sp \subset \Pn {15}$: Relations among these classes arise from
the tautological sequence, the conormal bundle sequence
$$ 0 \to {\cal N}^*_{\Sp\subset\Pn {15}} \to \Omega^1_{\Pn {15}}\mid_{\Sp} \to \Omega^1_{\Sp} \to 0$$ 
and ${\cal N}^*_{\Sp\subset\Pn {15}} \cong B^*(-2h)$.

With $b_i = c_i(B)$ theses relations give $b_1 = -2h, b_2=2h^2, b_4=-2h^4-2hb_3, b_5=0$ and 
$$b_3^2+8b_3h^3+8h^6=0=6h^5b_3+7h^8.$$
The last 2 equations define a complete intersection in $\Q[h,b_3]$. 
Since $$h^{10}=12 \in H^{20}(\Sp,\Q) \cong \Q,$$
there are no further relations between $h$ and $b_3$. 
To deduce
$$H^*(\Sp,\Q) \cong \Q[h,b_3]/(b_3^2+8b_3h^3+8h^6,6h^5b_3+7h^8)$$
we note that the odd cohomology of $\Sp$ vanishes and compute the Euler number of $\Sp$:
$$e(\Sp)=c_{10}(\Omega^1_{\Sp}) = {4 \over 3}h^{10} = 16,$$
which equals the length of the right hand side ring above. 
\bigskip

 $G \cong \Gr  (3,B)$. So $H^*(G,\Q)$ is generated as $H^*(\Sp,\Q)$-algebra by the
Chern classes $u_1,u_2,u_3$ of the universal subbundle U and relation follow from
the exact sequence
$$0 \to U \to B_G \to Q \to 0.$$
This gives
$$H^*(G,\Q) \cong H^*(\Sp,\Q)[u_1,u_2]/(f,g)$$
with
$$f=h^4-h^2u_2-{1 \over 2}u_2^2-{1 \over 2}b_3u_1+2h^3u_1-2hu_2u_1+3h^2u_1^2-{1 \over 2}u_2u_1^2
 +2hu_1^3+{1 \over 2}u_1^4$$
and 
$$g = b_3h^2-{1 \over 2}b_3u_2+2h^3u_2-hu_2^2+b_3hu_1+3h^2u_2u_1-2u_2^2u_1$$ 
   $$ -{1 \over 2}b_3u_1^2+2h^2u_1^3+{1 \over 2}u_2u_1^3+2hu_1^4+{1 \over 2}u_1^5.$$

\bigskip

$pr_2 \colon G \to M$ is a $\Pn 1$-bundle: $G \cong \P(E)$ with 
$E = pr_{2*}pr_1^*{\cal O}(h)$. Hence
$$H^*(G,\Q) \cong H^*(M,\Q)[h]/(h^2-c_1h+c_2)$$
with $ c_1=c_1(E),c_2=c_2(E)$. To compute these classes we note that $u_1,u_2,u_3$
are classes on $M$, since $U$ is induced from $M \subset \Gr  (3,10)$ and compute
$$Ann(u_1^{14}) = (u_1^2,u_2,h^2+u_1h) \subset H^*(G,\Q).$$
Since $H^4(M,\Q) \subset H^4(G,\Q)$ is 3- respectively 4-dimensional
$u_1^2,u_2,h^2+u_1h \in H^4(M,\Q)$. So 
$$c_2 = -h^2-u_1h, \quad c_1=-u_1.$$

\bigskip

For canonical classes we obtain
$$K_{\Sp} = c_1(\Omega^1_{\Sp}) = -8h, K_{G/\Sp} =c_1({\cal H}om(U,Q)^*)= 6h+5u_1,$$ 
hence $K_G=-2h+5u_1.$  $K_{G/M} = -2h +c_1(E)=-2h-u_1$, so
$$K_M = 6u_1.$$

\bigskip
$M_Y \subset M$ , the set of lines in $Y= \Sp \cap \Pn {10} \subset \Pn {15},$ is defined by the
vanishing of 5 sections of $E = pr_{2*}pr_1^*{\cal O}(h)$. Hence
$[M_Y]=(-h^2-u_1h)^5$ and $K_{M_Y}=K_M+5c_1(E)=u_1$. So
$M_Y$ is Fano with anti-canonical class $-u_1$ induced by the Pl\"ucker embedding
of $M_Y \subset M \subset \Gr  (3,10) \subset \Pn {119}$. Since
$$(-h^2-u_1h)^5(-u_1)^5h=11u_1^6h^{10} \in H^*(G,\Q),$$
the degree of the anti-canonical embedding of $M_Y$ is $11 \cdot 5 \cdot 12 =660$. 
In particular $-K_{M_Y}$ is not divisible in $Pic(M_Y)$, i. e. $M_Y$ is a Fano variety
of index 1.$\slutt$\medskip

Theorem (0.8) now follows from (8.6) and (9.2).

\vskip30pt

{\bf References}\bigskip

{\parindent 75pt

\item {[Alexander, Hirschowitz]}{\it Alexander, J., Hirschowitz, A.:}  Polynomial interpolation
in several variables, J. of Alg. Geom. {\bf 4} (1995), 201-222   
\item {[ACGH]}{\it Arbarello, E., Cornalba M. , Griffiths P. A., Harris J.:} Geometry of algebraic curves.
Vol. I. Grundlehren der Mathematischen  Wissenschaften, {\bf 267}, Springer-Verlag, New York-Berlin 1985
\item {[Buchweitz]} {\it Buchweitz, R. O.:}
Contributions a la theorie des singularities, thesis, Paris 7, 1981
\item {[Chevalley]}{\it  Chevalley, C.:} The algebraic theory of spinors. Columbia University
Press, New York, 1954
\item {[Clebsch]}{\it Clebsch, A.:}  \"Uber Curven vierter Ordnung, J. Reine Angew. Math. {\bf 59} (1861), 125-145
\item {[Dixon]}{\it Dixon A.:} The canonical forms of the ternary sextic and the quarternary
quartic, Proc. London Math. Soc. (2) {\bf 4} (1906), 223-227 
\item {[Dixon, Stuart]}{\it Dixon A., Stuart T.:} On the reduction of the ternary quintic and
septimic to their canonical forms, Proc. London Math. Soc. (2) {\bf 4} (1906), 160-168 \item
{[Dolgachev, Kanev]} {\it Dolgachev, I, Kanev V.:} Polar covariants of plane cubics and
quartics.  Adv. Math. {\bf 98}  (1993), 216-301 
\item {[Ein]}{\it Ein, L.:}  Varieties with
small dual varieties. I. Invent. Math. {\bf 86} (1986), 63-74 
 \item {[Eisenbud]}{\it Eisenbud, D.:} Commutative
algebra. With a view toward algebraic geometry. GTM {\bf 150} Springer-Verlag, New York, 1995
\item {[Eusen, Schreyer]}{\it Eusen F., Schreyer F.-O.:}  A remark to a conjecture of Paranjape
and Ramanan, Preprint, Bayreuth 1995
 \item {[Fano]}{\it Fano, G.:} Sulle variet\`a algebriche a tre dimensioni
a curve-sezioni canoniche,  Mem. R. Acc. d'Italia {\bf 8} (1937)  23-64 
\item {[Hilbert]}{\it Hilbert, D.:} 
Letter adresse\'e \`a M. Hermite, Gesam. Abh. vol. II, pp. 148-153 
\item {[Iarrobino]}{\it Iarrobino, A.:}
Compressed algebras: Artin algebras having socle degree and maximal length, Trans. Amer. Math. Soc.
{\bf 285} (1984), 337-378  
\item {}  {\it Iarrobino, A.:} 
Associated Graded Algebra of a Gorenstein Artin Algebra, Memoirs of AMS, {\bf 514}, (1994) 
\item {}{\it Iarrobino, A.:}  Inverse systems of a symbolic power. II:  The Waring problem
for forms, J. of. Algebra, {\bf 174}  (1995), 1091-1110 
 \item {[Iarrobino, Kanev]}{\it Iarrobino, A. Kanev, V.:} The length of a homogeneous
form, determinental loci of catalecticants and Gorenstein algebras, manuscript May 1996 
\item {[Kn\"orrer]}{\it Kn\"orrer H.:}  Cohen-Macaulay modules on hypersurface singularities I, Invent.
Math. {\bf 88} (1987), 153-164
\item {[Kustin, Miller]}{\it Kustin A. R., Miller M.:} Classification of the Tor-algebras of codimension
four Gorenstein local rings, Math. Z. {\bf 190} (1985), 341-355 
\item {[Lazarsfeld, Van de Ven]} {\it Lazarsfeld R., Van de Ven A.:}  Recent work of F.L.Zak,
DMV Seminar Band 4, Birkhauser (1984)  
\item {[Macaulay]}{\it Macaulay, F.S.:}   Algebraic theory of modular systems. Cambridge
University Press, London, (1916) 
\item {[MACAULAY]}{\it Bayer, D., Stillman, M.:}  MACAULAY: A system for computation in
algebraic geometry and commutative algebra, Source and object code available for Unix and Macintosh
computers.  Contact the authors, or download from zariski.harvard.edu via anonymous ftp.  
\item {[Mukai]}{\it Mukai, S.:} Curves $K3$ surfaces and Fano 3-folds of genus$\leq 10$, in "Algebraic
Geometry and Commuatative Algebra in Honor of Masayoshi Nagata", Kinokuniya,
Tokyo (1988), 357-377
\item {}{\it  Mukai, S.:}  Fano 3-folds,  Complex Projective Geometry, London Math. Soc. L. N. S. {\bf
179}, Cambridge University Press (1992), 255-263 
\item {}{\it  Mukai, S.:}  Polarized $K3$ surfaces of genus 18
and 20,  Complex Projective Geometry, London Math. Soc. L. N. S. {\bf 179}, Cambridge
University Press (1992), 264-276 
\item {}{\it  Mukai, S.:} Curves and Symmetric Spaces, I, American J. Math.  117,
(1995) 1627-1644 
\item {[Palatini]}{\it Palatini, F.:} Sulla rappresentazione delle forme ternarie
mediante la \break somma di potenze di forme lineari, Rom. Acc. L. Rend. (5) {\bf 12} (1903), 378-384 
\item {[Reye]}{\it Reye, T.:} Erweiterung der Polarenteori algebraische Fl\"achen, J. Reine
\break Angew. Math. {\bf 78} (1874), 97-114 
\item {}{\it Reye, T.:} Geometrischer Beweis des {\it Sylvester}schen Satzes: "Jede 
quater-\break n\"are cubische Form ist darstellbar als Summe von f\"unf Cuben linearer Formen", J.
Reine Angew. Math. {\bf 78} (1874), 114-122 
\item {}{\it Reye, T.:} Darstellung quatern\"arer biquadratischer
Formen als Summen von zehn Biquadraten, J. Reine Angew. Math. {\bf 78} (1874), 123-129 
\item {[Richmond]}{\it Richmond, H. W.:} On canonical forms,  Quart. J. Pure Appl. Math. {\bf
33} (1904), 331-340 
\item {[Room]} {\it Room, T. G.:} A synthesis of the Clifford matrices and its generalizations.
American J. Math. (1952),  967-984  
\item {[Rosanes]}{\it Rosanes J.:} \"Uber ein prinzip der Zuordnung
algebraischer Formaen, J. \"uber r. u. angew. Math. {\bf 76} (1873), 312-330 
\item {[Salmon]}{\it Salmon,
G.:} Modern Higher Algebra, 4. Edition.  Hodges, Figgis, and Co., Dublin (1885) 
\item {[Schreyer]}{\it Schreyer F.-O.:}  Syzygies of canonical curves and special linear
series.  Math. Ann. {\bf 275} (1986), 105-137 

\item {[Scorza]}{\it Scorza, G.:} Sopra la teoria delle figure polari delle curve piane del 4.
ordine, Ann. di Mat. (3) {\bf 2} (1899), 155-202 
\item {[Scorza]}{\it Scorza, G.:} Un nuovo theorema sopra
le quartiche piane generali, Math,. Ann. {\bf 52} (1899), 457-461 
\item {[Sylvester]}{\it Sylvester, J.J.:}  Sketch of a memoir on elimination, transformation and
canonical forms, Collected Works, Vol. I, Cambridge University Press, (1904), 184-197 
\item {} {\it Sylvester,
J.J.:}  An essay on canonical forms, supplemented by a sketch of a memoir on elimination,
transformation and canonical forms, Collected Works, Vol. I, Cambridge University Press,
(1904), 203-216 
\item {}{\it Sylvester, J.J.:} Sur une extension d'un th\'eor\`eme de Clebsch 
relatif aux courbes\break 
de quarti\`eme degr\'e, Compte Rendus de l'Acad. de Science {\bf 102} (1886), 1532-1534 (Collected
Math. Works IV p. 527-530) 
\item {[Terracini]}{\it Terracini, A.:} Sulle $V_k$ per cui la variet\`a degli
$S_h (h+1)$-seganti ha dimensione minore dell'ordinario.  Rend. Circ. Mat. Palermo {\bf 31} (1911),
392-396 
\item {[Zak]}{\it Zak, F. L.:} Varieties of small codimension arising from group actions, Addendum
to [Lazarsfeld, Van de Ven]

\item { }}

\bigskip
\vskip40pt
Authors' adresses:\bigskip
Kristian Ranestad\par
{\it Matematisk Institutt, UiO\par
 P.B. 1053 Blindern,\par
N-0316 Oslo, NORWAY}\par
e-mail: ranestad@math.uio.no\par
\bigskip
 Frank-Olaf Schreyer\par
{\it         
Universit\"at Bayreuth  \par
Mathematik          
D-95440 Bayreuth\par}
e-mail: schreyer@btm8x2.mat.uni-bayreuth.de

\bye